\RequirePackage{amsthm}

\documentclass[sn-mathphys-num,pdflatex]{sn-jnl}

\usepackage{graphicx}%
\usepackage{multirow}%
\usepackage{amsmath,amssymb,amsfonts}%
\usepackage{amsthm}%
\usepackage{mathrsfs}%
\usepackage[title]{appendix}%
\usepackage{xcolor}%
\usepackage{textcomp}%
\usepackage{manyfoot}%
\usepackage{booktabs}%
\usepackage{algorithm}%
\usepackage{algorithmicx}%
\usepackage{algpseudocode}%
\usepackage{listings}%
\usepackage{placeins}
\usepackage{paralist}
\usepackage{subfigure}

\theoremstyle{thmstyleone}%
\newtheorem{theorem}{Theorem}
\newtheorem{proposition}[theorem]{Proposition}%
\newtheorem{lemma}[theorem]{Lemma}%

\theoremstyle{thmstyletwo}%
\newtheorem{example}{Example}%
\newtheorem{conjecture}{Conjecture}%

\theoremstyle{thmstylethree}%

\newcommand{\R}{\mathbb{R}}   
\newcommand{\N}{\mathbb{N}}   
\newcommand{\dom}{\operatorname{dom}}
\newcommand{\co}{\operatorname{conv}}
\newcommand{\cone}{\operatorname{cone}}
\newcommand{\cartesian}{\ensuremath{\times}}
\newcommand{\eg}{\textit{e.g., }}
\newcommand{\ie}{\textit{i.e., }}

\usepackage{silence}
\begin{document}

\title[Conjugate of PLQ functions]{A linear-time algorithm to compute the conjugate of nonconvex bivariate piecewise linear-quadratic functions}


\author[1]{\fnm{Tanmaya} \sur{Karmarkar}}\email{tanmayak@student.ubc.ca}

\author*[1]{\fnm{Yves} \sur{Lucet}}\email{yves.lucet@ubc.ca}
\equalcont{These authors contributed equally to this work.}

\affil*[1]{\orgdiv{Computer Science, CMPS, I. K. Barber Faculty of Science}, \orgname{UBC Okanagan}, \orgaddress{\street{3187 University Way}, \city{Kelowna}, \postcode{V1T 1T7}, \state{BC}, \country{Canada}}}


\abstract{We propose the first linear-time algorithm to compute the conjugate of (nonconvex) bivariate piecewise linear-quadratic (PLQ) functions (bivariate quadratic functions defined on a polyhedral subdivision). Our algorithm starts with computing the convex envelope of each quadratic piece obtaining rational functions (quadratic over linear) defined over a polyhedral subdivision. Then we compute the conjugate of each resulting piece to obtain piecewise quadratic functions defined over a parabolic subdivision. Finally we compute the maximum of all those functions to obtain the conjugate as a piecewise quadratic function defined on a parabolic subdivision. The resulting algorithm runs in linear time if the initial subdivision is a triangulation (or has a uniform upper bound on the number of vertexes for each piece).

Our open-source implementation in MATLAB uses symbolic computation and rational numbers to avoid floating-point errors, and merges pieces as soon as possible to minimize computation time.}



\maketitle

\section{Introduction}\label{s:intro}

There are two motivations for the present work: obtaining tighter lower bounds for relaxation in global optimization, and computing the conjugate in computational convex analysis. Beyond applications, we want to understand the structure of the conjugate of piecewise linear-quadratic (PLQ) functions (bivariate functions defined on a union of polyhedral set on each of which the restriction of the function is quadratic). Understanding the structure is of interest in itself; it is also of interest to make progress toward computing the convex envelope (computing the conjugate is the first step toward computing the biconjugate, which is the convex envelope). The convex envelope provides the tightest relaxation when one wishes to solve a global optimization problem (the set of global minima of a function is included in the set of minima of its convex envelope). A large literature is available on global optimization and the search for tight relaxation; we refer to~\cite{LOCATELLI-13} for a detailed introduction and more details on relaxations relevant for our context.

Computing the convex envelope of a function is a hard problem; even computing the convex envelope of a multilinear function over a unit hypercube is NP-Hard~\cite{CRAMA-89}. However, results for specific functions exist in the literature, in particular for quadratic bivariate polynomials~\cite{LOCATELLI-14,LOCATELLI-14a,LOCATELLI-16,LOCATELLI-18a,KHADEMNIA-24}, and for convex envelopes of bilinear functions over triangles, rectangles and special polytopes~\cite{AL-KHAYYAL-83, SHERALI-90, LINDEROTH-05, ANSTREICHER-10, ANSTREICHER-12}. It is an active subject of research~\cite{LOCATELLI-24}.

PLQ functions play a significant role in variational analysis~\cite[Section 10E, Section 11D, p. 440]{ROCKAFELLAR-98a} due to the availability of calculus rules~\cite[Example 11.28, Proposition 11.32, Corollary 11.33, Proposition 12.30]{ROCKAFELLAR-98a}, and to their duality property~\cite[Theorem 11.42, Example 11.43, Theorem 11.42]{ROCKAFELLAR-98a}. Compared to piecewise linear functions, PLQ functions capture the curve of the original function more accurately
. The set of convex PLQ functions is closed under common convex operators, in particular under the Legendre-Fenchel transform; see~\cite[Page 484]{ROCKAFELLAR-98a},~\cite[Proposition 3, p. 17]{GOEBEL-00} and \cite{GARDINER-13}. 

The early idea of the computation of convex transforms can be traced back to~\cite{MOREAU-65}. However, development of most of the algorithms in computational convex analysis began with the computation of the conjugate with the Fast Legendre Transform (FLT)~\cite{BRENIER-89a}, which is studied in~\cite{CORRIAS-93a,LUCET-96a}. A linear-time algorithm is introduced in~\cite{LUCET-97b}. Those algorithms handle nonconvex functions but are restricted to grid domains. Extension to nongrid domains for convex functions only are proposed in~\cite{HAQUE-18,GARDINER-13} with a generalization to the partial conjugate in~\cite{GARDINER-13}. Computation of the conjugate of convex univariate PLQ functions has been well studied and linear time algorithms have been developed, both for the PLQ~\cite{LUCET-06} and the graph-matrix (GPH) models~\cite{GARDINER-11a}. All in all, there is no algorithm to compute the conjugate of nonconvex PLQ functions over nongrid domains.

Implementation of conjugate computation algorithms can be found in the CCA library~\cite{LUCET-17,LUCET-21} that contains efficient algorithms to manipulate convex functions and to compute fundamental convex analysis transforms arising in the field of convex analysis. Algorithms to compute the conjugate numerically on grids have  been based on either parameterization~\cite{HIRIART-URRUTY-06}, manipulation of graphs (GPH model)~\cite{GARDINER-13}, or the computation of the Moreau envelope~\cite{LUCET-05c}. More complex operators such as the proximal average operator~\cite{BAUSCHKE-07a,BAUSCHKE-06a} can be built by using a combination of addition, scalar multiplication, and conjugacy operations. 

A numerical library to determine in linear time whether a PLQ function is convex~\cite{SINGH-21} is available in MATLAB~\cite{CCA2}. This library is the first to handle general piecewise quadratic functions that are not necessarily continuous and are defined on any polyhedral subdivision without approximating them with a grid.

Computational Convex Analysis has many applications in the fields of image processing, network communication, PDE, geographic information systems, computer aided design, molecular biology, medical imaging, computer graphics and robotics; see the survey~\cite{LUCET-10}.

We propose the first algorithm to compute the conjugate of a general bivariate PLQ function not necessarily convex without approximating it with a grid. Our algorithm is split in three steps. In step 1, we apply \cite{LOCATELLI-16} to compute the convex envelope of each piece; in step 2, we apply \cite{KUMAR-19} to compute the conjugate of each resulting piece; and in step 3, we compute the maximum of all the resulting piecewise functions. 

Our contribution is fourfold: 
 \begin{inparaenum}[(1)]
    \item proving that the conjugate is a piecewise quadratic function defined on a parabolic subdivision (Kumar~\cite{KUMAR-19} had in addition a fractional form; we prove that such fractional form cannot occur);
    \item performing the very last step to obtain the formula for the conjugate (that step was not performed in~\cite{KUMAR-19}), and proving that computing the maximum results in a piecewise quadratic function defined on a parabolic subdivision; 
    \item proving the entire algorithm runs in linear time; and 
    \item releasing an open-source code that includes numerous techniques to minimize computation time (merging adjacent pieces with equal functions using symbolic computation).
\end{inparaenum}

We begin by giving some preliminaries in Section~\ref{s:prelim} and the overall algorithm in Section~\ref{s:alg}. We present examples in Section~\ref{s:examples}, and conclude with future work in Section~\ref{s:conclusion}.

\section{Preliminaries}\label{s:prelim}
The domain of a function $f:\R^n\to \R\cup\{+\infty\}$ is the set $\{x : f(x)<+\infty\}$, and the function is called proper if it has nonempty domain. The indicator function of $C \subset \R^n$ is denoted $I_C:\R \to (-\infty,\infty]$ with $I_C(x)=0$  if $x \in C$ and $+\infty$ otherwise. 

A function is called a piecewise function if it can be written $f(x)=f_i(x)$ for $x\in R_i$ with 
$\cup_i R_i=\R^n$. We call  the pair $(f_i,R_i)$ a piece. 
Recalling~\cite[Definition 10.20]{ROCKAFELLAR-98a}, a function $f:\R^n\to \R\cup\{+\infty\}$ is called Piecewise Linear-Quadratic (PLQ) if $\dom f = \cup_i R_i$ and $f$ restricted to $R_i$, noted $f_i$ can be written $f_i(x)=1/2 x^T Q x + q^T x + \kappa$ with $\kappa\in \R$, $q\in\R^n$, and $Q\in\R^{n\times n}$ a symmetric matrix. We use column vectors with $q^T x$ denoting the dot product of row vector $q^T$ with column vector $x$.
PLQ functions enjoy several properties, \eg their domain is closed, they are lower semicontinuous, and continuous relative to their domain~\cite[Proposition 10.21]{ROCKAFELLAR-98a}.

\begin{example}[Univariate PLQ function.]\label{uPLQ}
    The function $f(x) = x^2$ if $x \le 0$, $f(x) = 1-(1-x)^2$ if $0 < x < 1$, and       $f(x) = x^2$ when $x \ge 1$ is an example of a univariate PQ function.
\end{example}

\begin{example}[Bivariate PLQ function.]\label{bPLQ}
    The function 
    \begin{align*}
        f(x,y) &= 2x^2-xy-y^2,   &  (x,y) \in \co ((2,1),(3,5),(6,3)), \\
        f(x,y) &= x^2+xy-y^2,  &  (x,y) \in \co ((2,1),(6,3),(4,0)), \\
    \end{align*}
    is an example of a PLQ function in two variables.
\end{example}

An underestimator is any function that lies below a given function. The closed convex envelope is the largest lower semicontinuous convex underestimator of $f$.

\begin{example}[Convex Envelope of Example \ref{uPLQ}]
    The convex envelope $\co f$ of the function $f$ defined in Example \ref{uPLQ} is $\co f(x) = x^2$ if $x \le 0$, $\co f(x) = x$ if $0 < x < 1$, and $\co f(x) = x^2$ if $x \ge 1$.
\end{example}

We note 
\[f^*(y)=\sup_{x\in \R^n}(x^Ty - f(x))\] 
the (Legendre-Fenchel) conjugate of $f$.  The biconjugate is noted $f^{**}$, and we recall that when $f$ is proper, closed and convex, $f=f^{**}$, otherwise $f^{**}$ is the closed convex envelope of $f$~\cite[Theorem 11.1]{ROCKAFELLAR-98a}.    

To describe the domain of piecewise functions, we note a hyperplane $H \subset \R^d$ a set of the form $H= \{x \in \R^d: \alpha^Tx - \beta = 0\}$, and a halfspace $H \subset \R^d$ a set of the form $H= \{x \in \R^d: \alpha^T x - \beta \le 0\}$ or $H= \{x \in R^d: \alpha^Tx - \beta \ge 0\}$; where the vector $\alpha \in \R^d \backslash  \{0\}$ is called a normal vector to $H$.

A polytope is a convex combination of finite numbers of vertices $\{v_1, \dots, v_k\}$. We write $
\co V = \{ \lambda_1 v_1 + \dots + \lambda_k v_k: \lambda_i \ge 0, \sum_{i=1}^{k}\lambda_i=1 \}$ or in matrix notation $\co V =  \{ V\lambda: \lambda_i \ge 0, \sum_{i=1}^{k}\lambda_i=1 \}$,
where $V$ is the matrix formed by column vectors. Note that by our definition, polytopes are always bounded.

A polyhedral cone is defined as a conic combination of a finite numbers of directions $\{d_1,\dots,d_m\}$; we write $\cone C = \{\mu_1d_1+ \dots + \mu_md_m : \mu_i \ge 0\}$, and $\cone C = \{D\mu : \mu_i \ge 0\}$;
where $D$ is the matrix formed by column vectors $d_i$.

A polyhedral set (polyhedron) $P \subset \R^d$ is defined as the intersection of a finite number of closed halfspaces and hyperplanes; it has the $H$-representation 
\[P=\{x\in \R^d: Ax\leq b\}\] 
with $A\in\R^{m\times d}$ and $b\in \R^m$. 
We can also describe a polyhedral set in Vertex-representation. According to the Minkowski-Weyl Theorem~\cite[Corollary 3.53]{ROCKAFELLAR-98a}, a polyhedral set can be written as a sum of a polytope and a polyhedral cone. Therefore, in $V$-representation a polyhedral set $P$ is described as
\[
	P = \co V + \cone C 
	  = \bigg\{ V\lambda + D\mu: \lambda \in \R^k, \lambda \ge 0, \mu \ge 0, \sum_{i=1}^{k}\lambda_i=1 \bigg\}.
\]
Since the intersection of convex sets is convex, polyhedral sets are convex.

A $d$-dimensional face of a convex set $C\subset \R^n$ is a $d$-dimensional convex subset $C'$ of $C$ such that every (closed) line segment in C with a relative interior point in $C'$ has both endpoints in $C'$~\cite[Page 162]{ROCKAFELLAR-70}. The empty set and $C$ itself are faces of $C$. The two-dimensional faces of $C$ are called faces.
The one-dimensional faces of a convex set C are called edges. For $u,v\in\R^d, u \neq v$ and $\delta=(\delta_1,\delta_2)\in \Delta\subset\R^2,$ an edge can be written as
\[\bigg\{x\in\R^d :   x=\delta_1u+\delta_2v, \delta_1+\delta_2=1\bigg\}.\]
An edge in $\R^2$ is a segment if $\Delta = \R_+ \cartesian \R_+ $ where $\R_+= \{\alpha\in \R:\alpha\ge 0\}$, a ray if $\Delta = \R_+ \cartesian \R$, and a line if $\Delta = \R^2$.
The zero-dimensional faces of a convex set $C$ are called vertices and are the extreme points of $C$.

A convex set defined as the union of a finite number of polyhedral regions, namely
$R = \bigcup_{i=1}^n R_{i}$,
where $R \subseteq \R^2$, is said to be a polyhedral subdivision~\cite[Definition 1]{PATRINOS-11}, if $R_i$ is a polyhedral set and for any  \(j, k \in \{1, \dots, n\}, j \neq k, R_{j} \cap R_{k}\) is either empty, a vertex or an edge.
(Polyhedral subdivisions eliminate degenerate subdivisions for which the intersection of two polyhedral regions is a strict subset of an edge, or of a face.)

Assume that $f: \R^2 \rightarrow \R \cup \{+\infty\}$  is a PLQ function and that $\dom f = \cup_i P_i $ is a polyhedral subdivision. An entity is a $d$-dimensional face of $P_i$ for some index $i$. An entity is either a vertex, an edge or a face of the polyhedral subdivision of $\dom f$~\cite{GARDINER-11}.

We now define the geometric shapes required to describe the domain of the conjugate functions we obtain. 
A parabola $P\subset\R^2$ is a subset of the plane that can be written as
\[P=\{(x,y)\in\R^2: ax^2+b x y+c y^2+ d x+e y +f = 0\},\]
where $a, b, c, d, e, f\in\R$ are not all zero and satisfy \(b^2 - 4ac = 0\) .
Note that our definition includes lines when $a=b=c=0$, and the empty set when $a=b=c=d=e=0$, $f\neq 0$; but excludes the entire plane since we impose that $(a,b,c,d,e,f)\neq 0$.

A parabolic region \(P_r \subset  \R^2\) is formed by the intersection of a finite number of parabolic inequalities. It can be written as 
\[Pr = \{x \in \R^2: a_ix^2 + b_ixy + c_iy^2 + d_ix + e_iy + f_i \leq 0, i=1, \dots , k\},\] 
where $(a_i,b_i,c_i,d_i,e_i,f_i)\neq 0$, and \( b_i^2- 4a_i c_i = 0\).
The following sets are special cases of parabolic regions: polyhedral sets, polyhedral cones, polytopes, edges, and vertices. The set $\{(x,y)\in\R^2: y\geq x^2\}$ is an example of a non-polyhedral parabolic region while the set $\{(x,y) : y\leq x^2, y\geq 1\}$ is a nonconvex non-connected parabolic region.

A convex set $R\subset \R^2$ is called a parabolic subdivision if $R$ can be written as the finite union of parabolic regions $R_i$, \ie $R = \cup_{i=1}^n R_{i}$ and for any $j, k \in \{1, \dots , n\}$, $j\neq k$, the intersection $R_j \cap R_k$ is either empty or is contained in a parabola.
A polyhedral subdivision is a special case of a parabolic subdivision (since lines are a special case of parabolas).

A face of a parabolic region \(P_r \subset  \R^2\) is defined as the interior of a nonempty intersection of a finite number of parabolic inequalities. It can be written as 
\[P_r = \{x \in \R^2: a_ix^2 + b_ixy + c_iy^2 + d_ix + e_iy + f_i < 0, i=1, \dots , k\},\] 
where $(a_i,b_i,c_i,d_i,e_i,f_i)\neq 0$, and \( b_i^2- 4a_i c_i = 0\).

Finally, we use standard definitions for subdifferentials, subgradients, and normal cones~\cite{BAUSCHKE-24,HIRIART-URRUTY-93a,HIRIART-URRUTY-93b,ROCKAFELLAR-70a}.

\section{Computing Conjugates}\label{s:computing_conjugate}
While the convex envelope of a piecewise function is not always the same as the convex envelope of each piece, for bivariate quadratic polynomials we can leverage the convex envelope of each piece to obtain the convex envelope. Let us denote the $i^{th}$ piece of a piecewise function $f$ by $(f_i,P_i)$ where $f_i$ is a function and $P_i$ is a set.
From \cite[Theorem 3.4.1]{HIRIART-URRUTY-93b}, we have $(\inf_i f_i)^*=\sup_i f_i^*$. Using \cite[Proposition 2.6.1]{HIRIART-URRUTY-93b} and the biconjugate theorem $\co f = f^{**}$, we get
\begin{align*}
	\co[\inf_i(f_i+I_{P_i})]=&\co[\inf_i(\co(f_i+I_{P_i}))],\\
	=&[\inf_i(\co(f_i+I_{P_i}))]^{**},\\
	=&[\sup_i([\co(f_i+I_{P_i})]^{*})]^*.
\end{align*}
Taking the conjugate on both sides, we obtain
\[f^* = \left( 	\co[\inf_i(f_i+I_{P_i})] \right)^* = \sup_i([\co(f_i+I_{P_i})]^{*}).\]

Hence we can find the conjugate of a bivariate PLQ function as follows.
\begin{enumerate}
	\item Compute the convex envelope of each piece $\co(f_i+I_{P_i})$.
	\item Compute the conjugate of each convex piece $[\co(f_i+I_{P_i})]^{*}$.
	\item Compute the maximum of the conjugates over the entire PLQ function to obtain
	$f^*=\sup\limits_{i}(\co(f_i+I_{P_i})^{*})$.
\end{enumerate}

Step 1 of this process, finding the convex envelope of each piece $(f,P)$ with $f(x)=x^TAx+b^Tx+c$, is computed as follow. We compute the eigenvalues of $A$. If they are both nonnegative, the function is convex and we are done. If they are both nonpositive, the function is concave and the convex envelope is obtained as the convex hull of the points $\{(x,f(x)) : x$ is a vertex of $P\}$. If $A$ is indefinite (one positive and one negative eigenvalue), we note that $\co f(x) = \co (x^TAx) + b^Tx+c$; see~\cite[P. 93]{HIRIART-URRUTY-93a}. Hence, we only focus on the quadratic part. We rotate the polyhedral set so that in the new basis, the function is now $x\mapsto x_1 \cdot x_2$. We then use~\cite{LOCATELLI-16} where a method to compute the convex envelopes of bilinear forms over general polytopes is given. 

Kumar~\cite{KUMAR-19a} showed how to compute the conjugate of a rational function defined on a polytope, which is Step 2 of the above process. 

Our focus is to explain how to perform Step 3, the last step required to obtain the conjugate, by computing the maximum of the functions obtained in Step 2 thereby obtaining the conjugate of a bivariate PLQ function. Repeating steps 2 and 3 to compute the biconjugate $f^{**}=(\sup\limits_{i}(\co(f_i+I_{P_i})^{*}))^*$ is left for future work.

We point out a significant simplification of the structure of the conjugate. Locatelli~\cite{LOCATELLI-16} showed that the convex envelope of a quadratic over a polytope is a piecewise function defined on a polyhedral subdivision where each function is either linear, quadratic, or a rational function (ratio of quadratic over linear). It turns out that these rational functions have a specific structure. Kumar~\cite{KUMAR-19a} subsequently showed that the conjugate of any of those functions is a piecewise function defined over a parabolic subdivision where the restriction to each piece is either linear, quadratic, or a fractional form 
\[
		g_f(s_1,s_2) = \frac{\psi_1(s_1,s_2)}{\zeta_{00}\sqrt{\psi_{\frac{1}{2}}(s_1,s_2)}} + \psi_0(s_1,s_2),
\]
where $\psi_1$, $\psi_0$ $\psi_{\frac{1}{2}}$ are linear functions in $s=(s_1,s_2)\in\R^2$ (justifications are in~\cite{KUMAR-19} and details in~\cite{KARMARKAR-24}). It turns out that the special structure of the rational functions considered prevent this general fractional form to occur.


\begin{proposition}
    Assume $f(x)=x^TAx+b^Tx+c +I_P$ with $P=\co\{x_1,x_2,x_3\}$ a triangular set in $\R^2$. Then its convex envelope $\co f = \min f_i + I_{P_i}$ is a piecewise function over a polyhedral subdivision where $f_i$ is quadratic, linear, or a rational function while $P$ is polyhedral. The conjugate of each piece of $\co f$ is a piecewise quadratic function defined over a parabolic subdivision.
\end{proposition}

\begin{proof}
Locatelli~\cite{LOCATELLI-16} showed that $\co f$ is a piecewise function over a polyhedral subdivision that can be linear, quadratic, or rational. Kumar~\cite{KUMAR-19} showed that the conjugate of each piece is either a quadratic function, or a fractional function defined over a parabolic subdivision. We only need to prove that the fractional form never occurs.

The fractional form only appears in the conjugate when we compute the conjugate of a rational function. We only get a rational function as the convex envelope when we solve the optimization problem from Locatelli with $\eta_h = -{(a+mb-q)^2}/{4m}-bq$ and $\eta_w$ linear. In this case $\eta_h$ corresponds to the convex edge $y=mx+c$ and $\eta_w$ corresponds to the vertex $(x_1,y_1)$. The optimization problem we are solving is
\[
    \max_{a,b} \{ \eta_w -ax - by \; : \;  \eta_w=\eta_h, (a,b) \in S_{r} \}
\]
which gives a solution which is a rational function
\[
r(x,y)=\frac{ax^2+bxy+cy^2+dx+ey+f}{gx+hy+k}
\]
where $a=-my_1$, $b=q$, $c=x_1$, $d=-qy_1+mx_1y_1$, $e=-qx_1-x_1y_1$, $f=qx_1y_1$, $g=-m$, $h=1$, and $k=-y_1 mx_1$. (The rational function is extended by continuity at vertex $(x_1,y_1)$.) The computation is verified symbolically in Appendix~\ref{a:nofractional} (file vertexNan.m).

From \cite{KUMAR-19}, the conjugate of the rational function of this form is a piecewise quadratic function defined on a parabolic division (a quadratic function is obtained for the conjugate corresponding to the convex edge while all other functions are linear). All other functions from the convex envelope, linear or quadratic, result in linear or quadratic functions for their conjugates. 
\end{proof}

Our implementation in MATLAB uses symbolic computation and rational numbers to avoid any floating-point errors. Using floating-point arithmetic can give rise to partitions of the domain with degenerate subsets, and floating-point coefficients of functions make it difficult to detect equal functions on adjacent domains. Indeed, one of the motivations of~\cite{SINGH-21} is to check that intermediate computations of conjugates are still convex (as they should be without floating point errors). After observing these difficulties, we decided to work with symbolic representation of our functions with rational coefficients.

\subsection{Algorithm}\label{s:alg}

We use the symbols $q$ for quadratic functions in the primal space, $r$ for rational functions in the primal space, and $s$ for quadratic functions in the dual space. We now fix our notations for each intermediate computation by naming each of the piecewise function involved.

The original PLQ function and its domain are noted $f$ and $\dom f$  where 
\[f = \min_{i=1 \dots n_q} (q_i + I_{\dom q_i}),\]
where $q_i$ is the quadratic function defining the $i^{th}$ piece, and $f$ has $n_q$ pieces. 

The convex envelope of $ q_i + I_{\dom q_i} $ is denoted
\[ \co (q_i + I_{\dom q_i}) = \min_{j=1 \dots n_{r_i}} (r_{i,j} + I_{\dom r_{i,j}}).\]
Here $r_{i,j}$ represents a rational function that defines the $j^{th}$ piece of the convex envelope of the $i^{th}$ piece of $f$.

The conjugate of $r_{i,j} + I_{\dom r_{i,j}}$ is written
\[(r_{i,j} + I_{\dom r_{i,j}})^* = \min _{k=1,\dots n_{s_{i,j}}} (s_{i,j,k} + I_{\dom s_{i,j,k}} ),\]
where $s_{i,j,k}$ represents a quadratic function which is the $k^{th}$ piece of $(r_{i,j} + I_{\dom r_{i,j}})^*$ .

The conjugate of $\co(q_i + I_{\dom q_i})$ is 
\begin{align*}
	[\co(q_i + I_{\dom q_i})]^*   & =\max_{j=1 \dots n_{r_i}} \min_{k=1,\dots n_{s_{i,j}}} (s_{i,j,k} + I_{\dom s_{i,j,k}} ) ,\\
	&= \min_{k=1 \dots n_{s_i}} (s_{i,k} + I_{\dom s_{i,k}} ).
\end{align*}
The last equality defines our notation. The conjugate is a piecewise function with $n_{s_i}$ pieces; each quadratic function is denoted $s_{i,k}$.

Finally, the conjugate of $f$ is denoted
\[f^* =  \max_{i=1 \dots n_q}  (s_{i,k} + I_{\dom s_{i,k}} ) =  \min_{k=1 \dots n_s} (s_k + I_{\dom s_k}). \]
It is a piecewise function with $n_s$ pieces where each quadratic function is denoted $s_k$.

While we use the same notation $s$, the number of indexes indicate which function it is referring: 1 index for a piece of $f^*$, 2 indexes for a piece of $[\co(q_i + I_{\dom q_i})]^*$, and 3 indexes for a piece of $(r_{i,j} + I_{\dom r_{i,j}})^*$.

\begin{figure}
	\centering
    \includegraphics[height=220px]{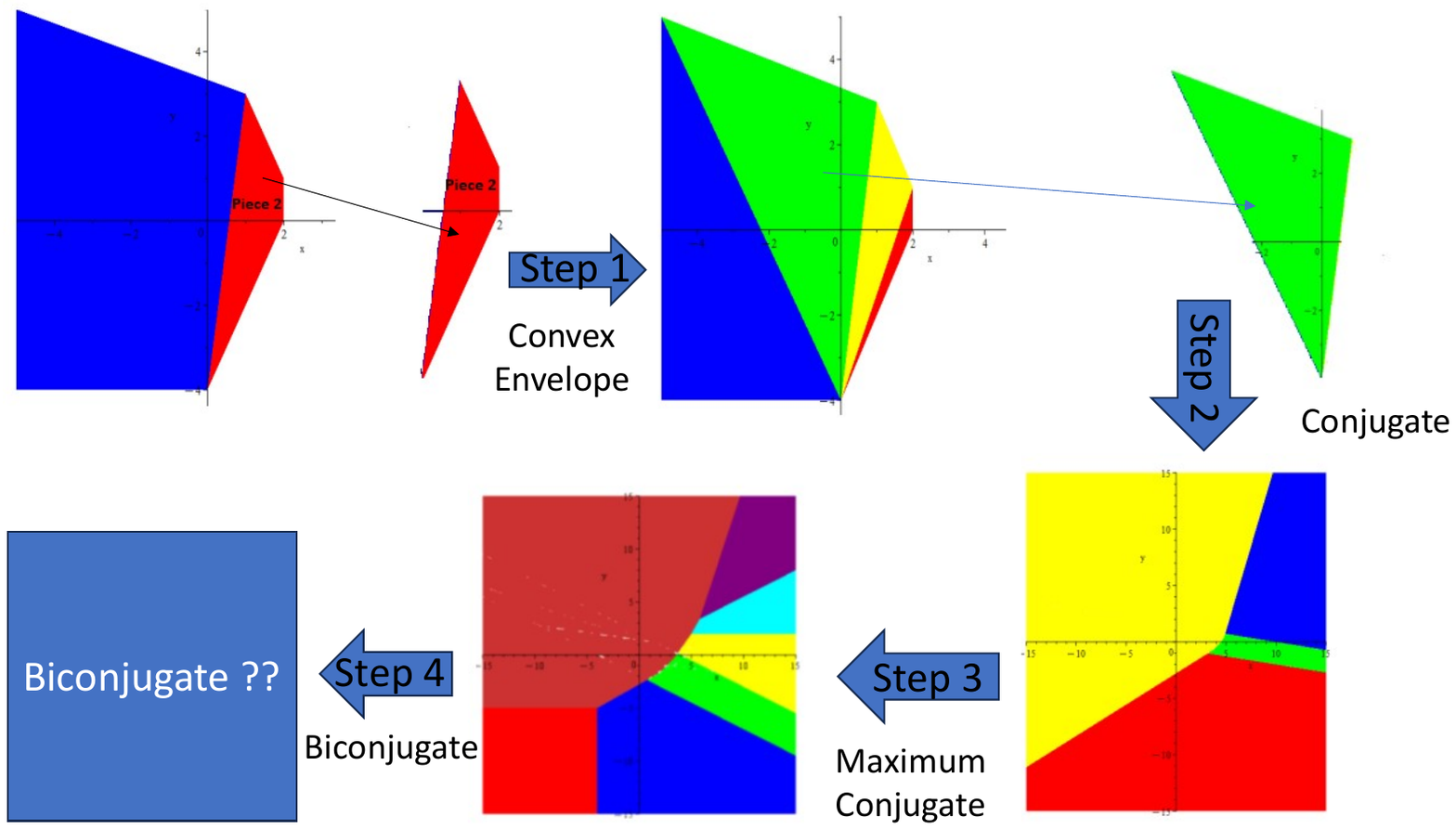}
	\caption{Illustration of the steps taken on the domains to compute the biconjugate of a PLQ function. Note that the biconjugate has a polyhedral subdivision, but an unknown explicit formula leading us to label the last box ``Biconjugate ??''.}
	\centering
	\label{fig:1}
\end{figure}

\begin{figure}
	\centering
	\subfigure[$ q_{i}.$]{
		\centering
        \includegraphics[height=70px]{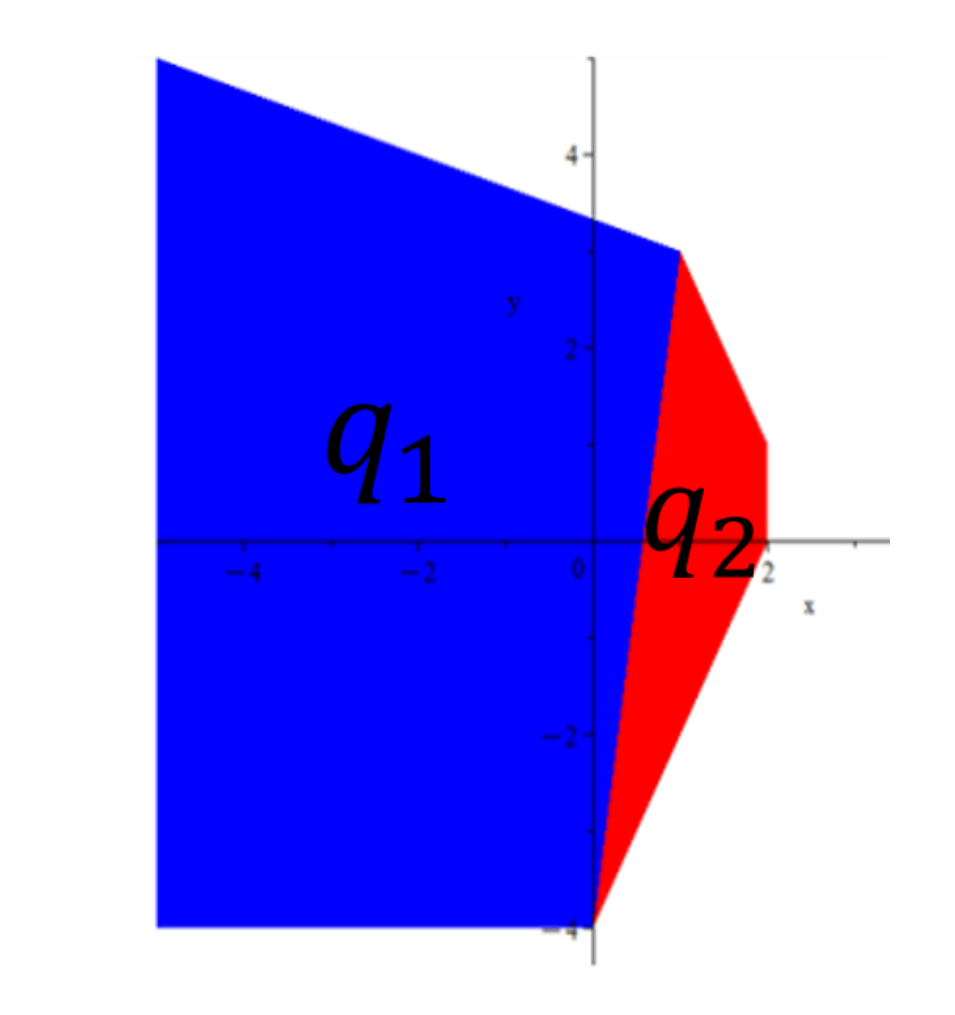}\label{f:f}}\\
	\subfigure[$r_{1,j}.$ ]{
		\centering
        \includegraphics[height=80px]{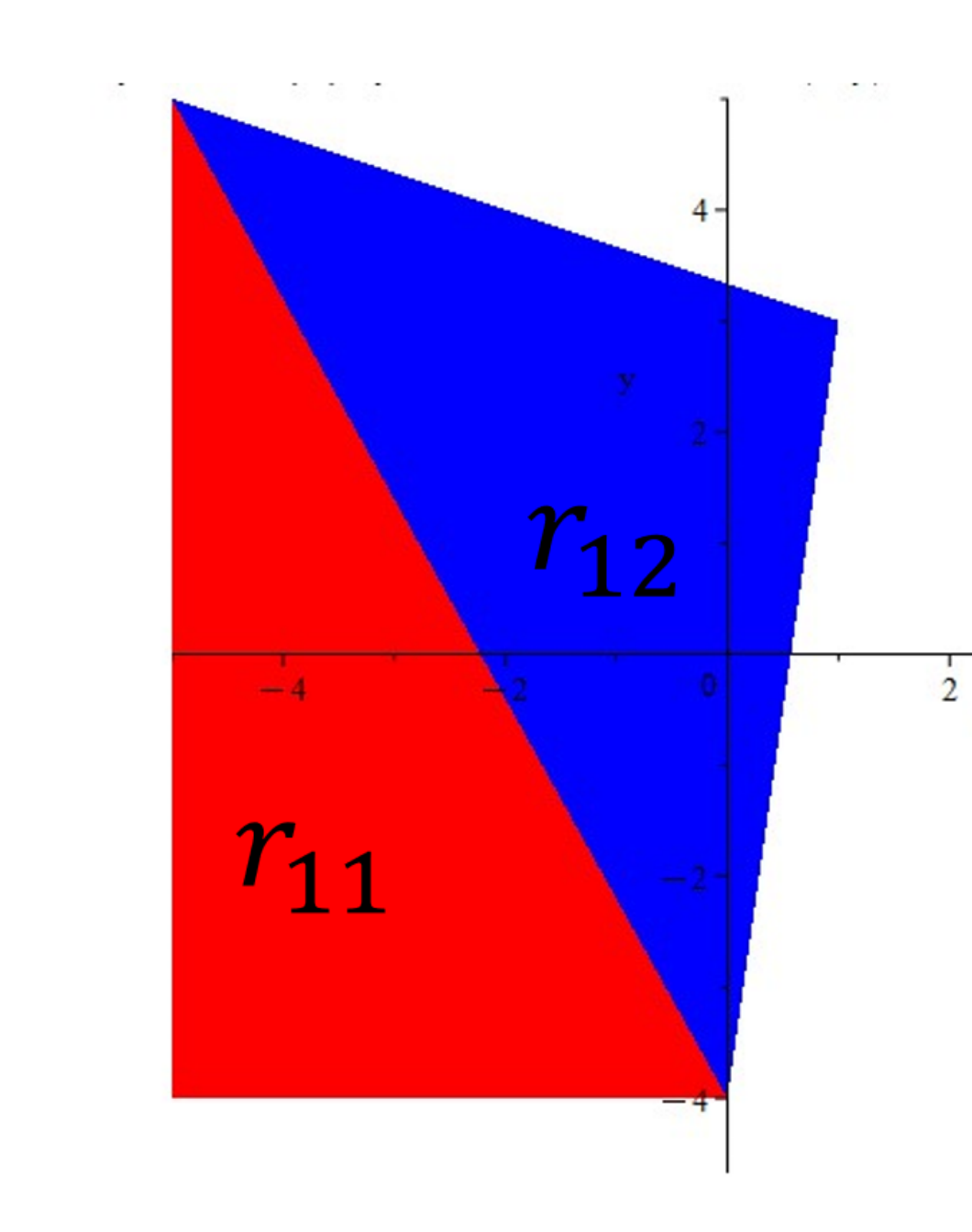}\label{f:f11}}
	\subfigure[$r_{2,j}.$ ]{
		\centering
        \includegraphics[height=80px]{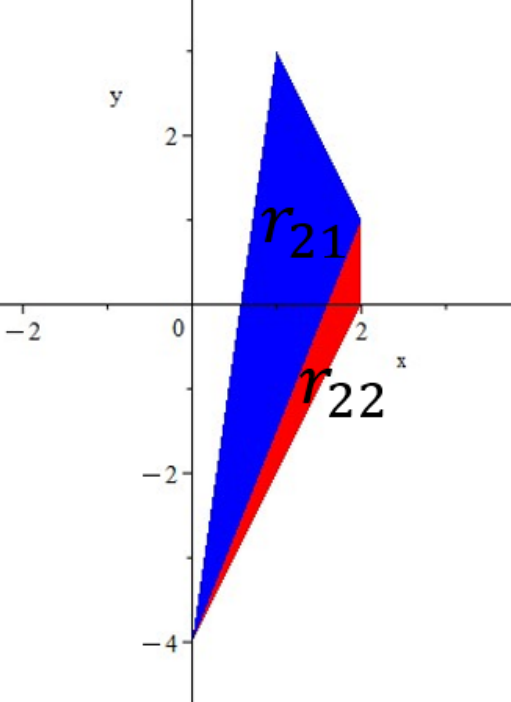}\label{f:f12}}\\
	\subfigure[$s_{1,1,k}.$ ]{
		\centering
        \includegraphics[height=80px]{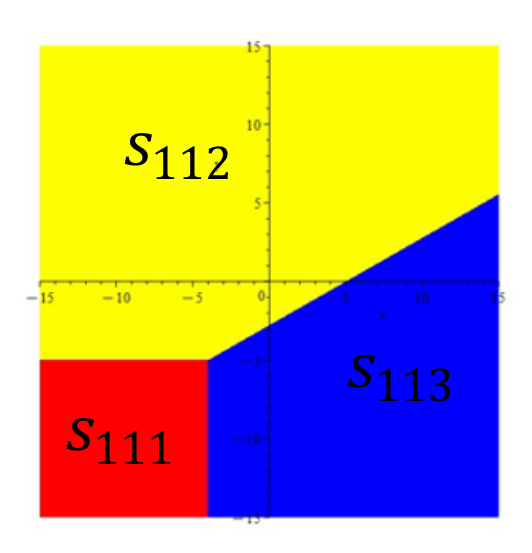}\label{f:conj12}}
	\subfigure[$s_{1,2,k}.$ ]{
		\centering
        \includegraphics[height=80px]{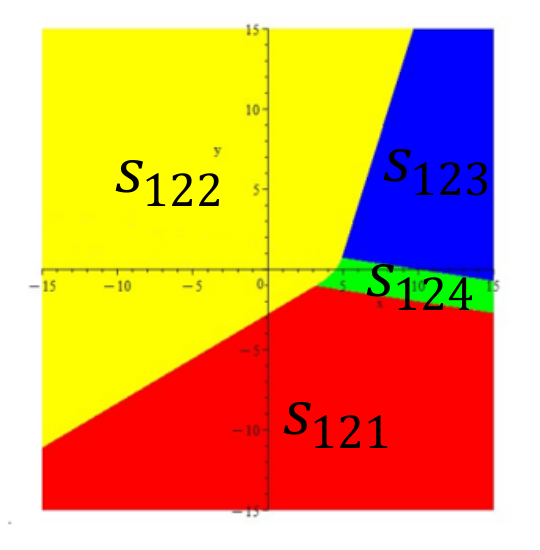}\label{f:conj11}}
	\subfigure[$s_{2,1,k}.$ ]{
		\centering
        \includegraphics[height=80px]{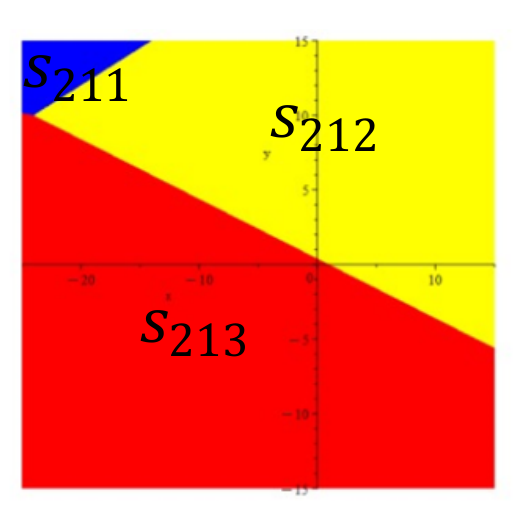}\label{f:conj22}}
	\subfigure[$s_{2,2,k}.$ ]{
		\centering
        \includegraphics[height=80px]{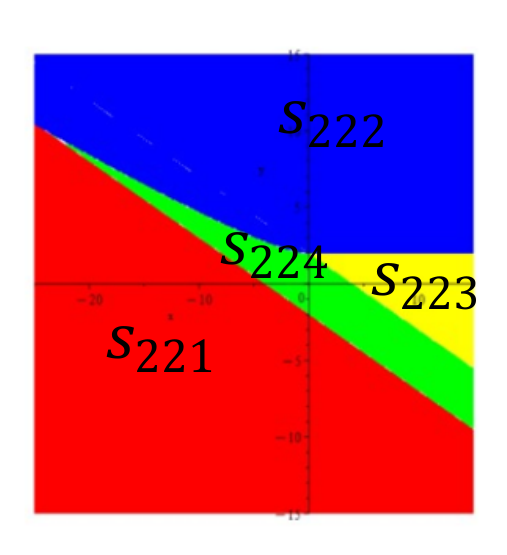}\label{f:conj21}}\\
	\subfigure[$s_{1,k}.$ ]{
		\centering
        \includegraphics[height=80px]{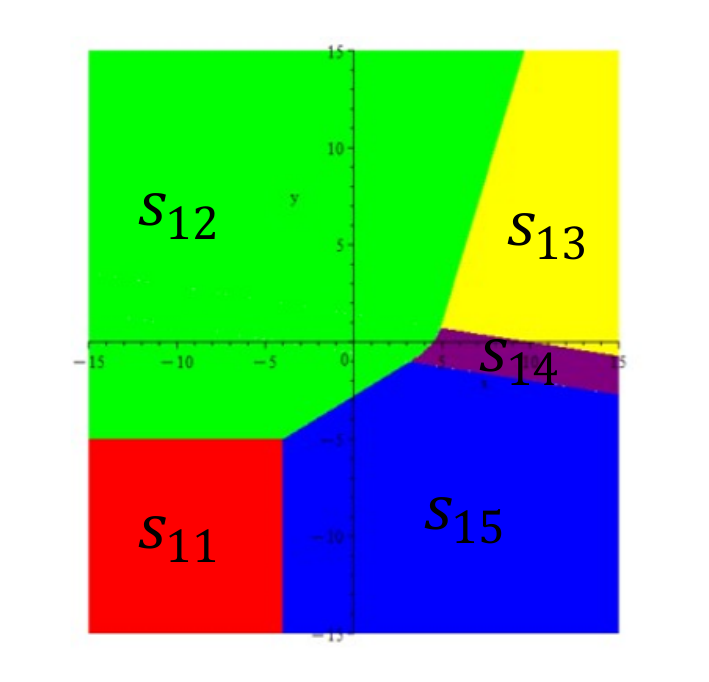}\label{f:conjmax1}}
	\subfigure[$s_{2,k}.$ ]{
		\centering
        \includegraphics[height=80px]{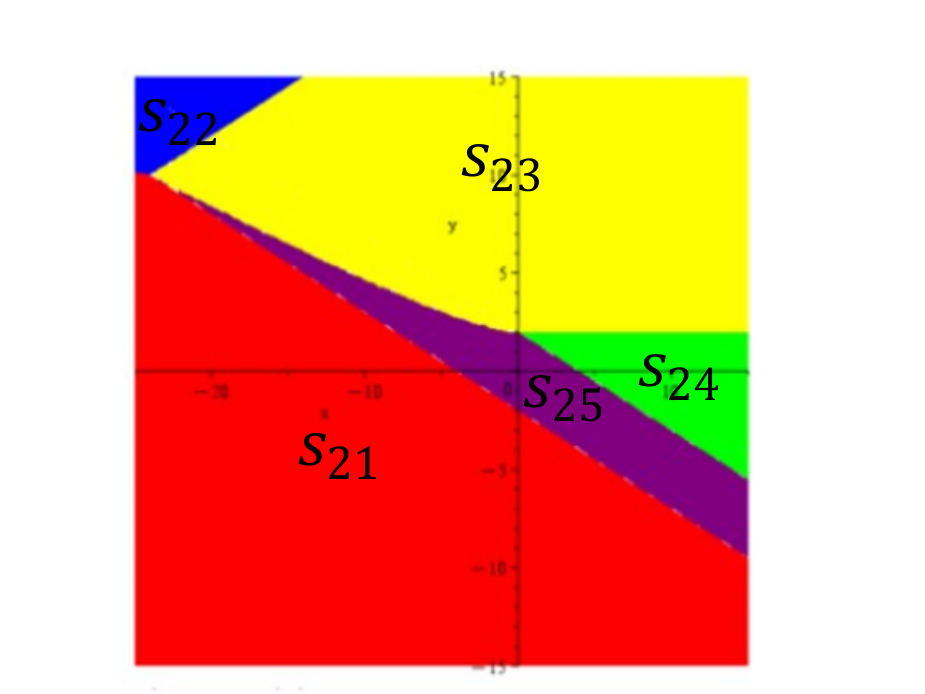}\label{f:conjmax2}}\\
	\subfigure[$s_{k}.$ ]{
		\centering
        \includegraphics[height=80px]{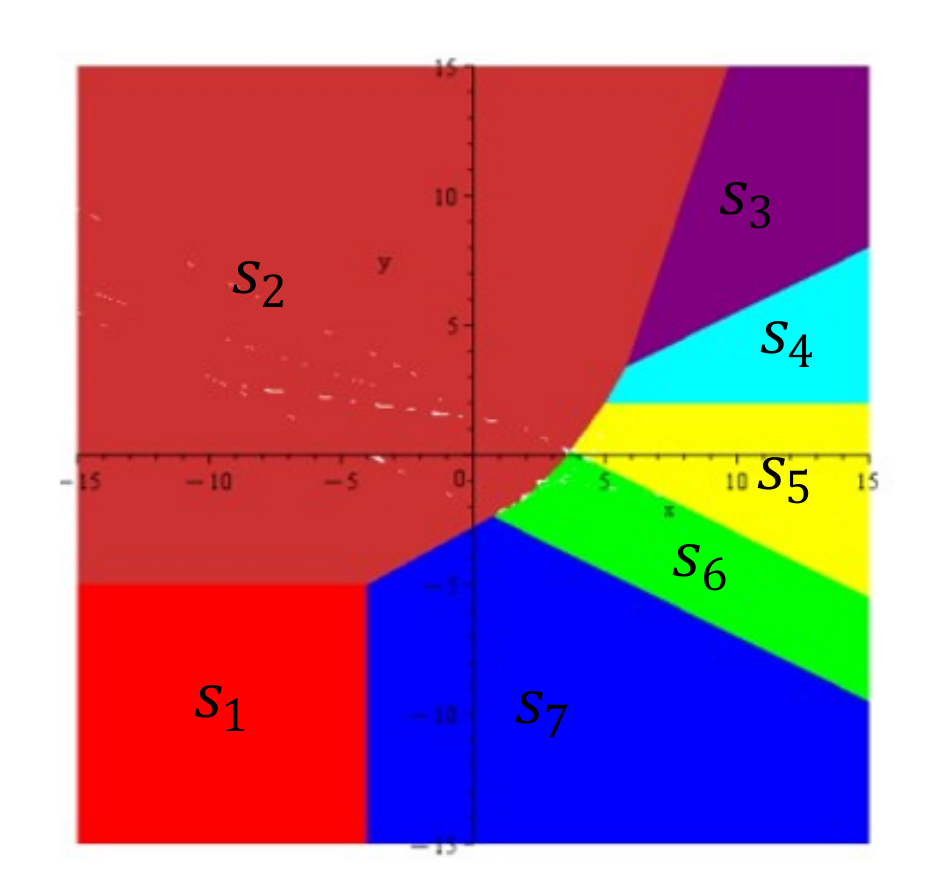}\label{f:conj}}
	\caption{Steps to compute the conjugate of a PLQ Function illustrated on the domain of the respective functions.(Functions defined in Table \ref{table:allsteps})}
	\label{fig:example1-all}
\end{figure}

We illustrate the overall flow of computing the conjugate of the PLQ function in figures~\ref{fig:1} and~\ref{fig:example1-all}. We are given a PLQ function with $n_q=2$ pieces. Each piece has a nonconvex bivariate function defined over a polyhedral region. As an example in Figure~\ref{fig:example1-all}, we have a PLQ function with two pieces. 

\begin{table}[tbph]
\centering
\caption{Function definitions for Figure \ref{fig:example1-all}.}
	\label{table:allsteps}
\begin{tabular}{@{}llr@{}} \toprule 
    Function name & Expression \\
    \midrule
	$q_1$ &$xy$ \\
    $q_2$ &$xy$  \\ 
    $r_{11}$ &$-4x-5y-20$  \\ 
    $r_{12}$ &$(155x-5y+4xy+35x^2+5y^2-100)/(7x-y+40)$  \\ 
    $r_{21}$ &$-24x + 10y + 40$  \\ 
    $r_{22}$ &$(8x+6y-4xy-2x^2+2y^2-8)/(y-2x+3)$  \\ 
    $s_{111}$ &$-5s_1+4s_2-20$  \\ 
    $s_{112}$ &$-5s_1+5s_2+25$  \\ 
    $s_{113}$ &$-4s_2$  \\ 
    $s_{121}$ &$-4s_2$  \\ 
    $s_{122}$ &$-5s_1+5s_2+25$  \\ 
    $s_{123}$ &$s_1+3s_2-3$  \\ 
    $s_{124}$ &$\frac{1}{28}s_1^2+\frac{1}{2} s_1 s_2+\frac{2}{7}s_1+\frac{7}{4}s_2^2-2s_2 +\frac{4}{7}$  \\ 
    $s_{211}$ &$s_1+3s_2-46$  \\ 
    $s_{212}$ &$2s_1+s_2-2$  \\ 
    $s_{213}$ &$-4s_2$  \\ 
    $s_{221}$ &$-4s_2$  \\ 
    $s_{222}$ &$2s_1+s_2-2$  \\ 
    $s_{223}$ &$2s_1$  \\ 
    $s_{224}$ &$\frac{1}{8}s_1^2+\frac{1}{2}s_1s_2+s_1+\frac{1}{8}1s_2^2-2s_2+2$  \\ 
    $s_{11}$ &$-5s_1+4s_2-20$  \\ 
    $s_{12}$ &$-5s_1+5s_2+25$  \\ 
    $s_{13}$ &$s_1+3s_2-3$  \\ 
    $s_{14}$ &$\frac{1}{28}s_1^2+\frac{1}{2} s_1 s_2+\frac{2}{7}s_1+\frac{7}{4}s_2^2-2s_2 +\frac{4}{7}$  \\ 
    $s_{15}$ &$-4s_2$  \\ 
    $s_{21}$ &$-4s_2$  \\ 
    $s_{22}$ &$s_1+3s_2-46$  \\ 
    $s_{23}$ &$2s_1+s_2-2$  \\ 
    $s_{24}$ &$2s_1$  \\ 
    $s_{25}$ &$\frac{1}{8}s_1^2+\frac{1}{2}s_1s_2+s_1+\frac{1}{8}1s_2^2-2s_2+2$  \\ 
    $s_{1}$ &$-5s_1+4s_2-20$  \\ 
    $s_{2}$ &$-5s_1+5s_2+25$  \\ 
    $s_{3}$ &$s_1+3s_2-3$  \\ 
    $s_{4}$ &$2s_1+s_2-2$  \\ 
    $s_{5}$ &$2s_1$  \\ 
    $s_{6}$ &$\frac{1}{8}s_1^2+\frac{1}{2}s_1s_2+s_1+\frac{1}{8}1s_2^2-2s_2+2$  \\ 
    $s_{7}$ &$-4s_2$  \\ 
	\bottomrule 
\end{tabular}
\end{table}
The algorithm follows the following steps
\begin{enumerate}
	\item Given $q_i$ compute the convex envelope of $q_i + I_{\dom q_i}$ to obtain $r_{i,j}$.
	\item From $r_{i,j}$, compute the conjugate of $r_{i,j} + I_{\dom r_{i,j}}$ to obtain $s_{i,j,k}$.
	\item From $s_{i,j,k}$, compute the conjugate of $\co(q_i + I_{\dom q_i})$ to obtain $s_{i,k}$, and deduce $s_k$ to obtain the conjugate of $f$.
\end{enumerate}

We apply Step 1 and get the convex envelope for each piece, as shown in Figure~\ref{f:f11} and Figure~\ref{f:f12}. This is followed by Step 2 to obtain the conjugate of each convex piece. The conjugates are shown in figures~\ref{f:conj12}-\ref{f:conj21} and finally the maximum of the conjugates are computed in Step 3. Figures~\ref{f:conjmax1}, \ref{f:conjmax2} show the conjugate of each piece and finally Figure~\ref{f:conj} shows the conjugate of the entire PLQ function.

\FloatBarrier 

From \cite[Theorem 1.1]{LOCATELLI-16}, we know that the convex envelope of $q_i + I_{\dom q_i}$ has a polyhedral subdivision where each piece is associated with a function of the form
\[
\frac{\alpha_6 x^2 + \alpha_5 y^2 + \alpha_4 xy + \alpha_3 x + \alpha_2 y + \alpha_1 }{\beta_3x + \beta_2y+\beta_1}.
\]
We refer to \cite{LOCATELLI-16} for the justification and to~\cite{KARMARKAR-24} for implementation details using our notations.

Step 2 computes the conjugate of each rational function over a polytope $r_{i,j} + I_{\dom r_{i,j}}$ to obtain a piecewise function defined on a parabolic subdivision and one of the following functional forms
\begin{equation}\label{eq:c2}
	\begin{aligned}
		g_q(s_1,s_2) &= \zeta_{11}s_1^2 + \zeta_{12}s_1s_2 + \zeta_{22}s_2^2 + \zeta_{10}s_1+\zeta_{01}s_2 + \zeta_{00},\\
		g_l(s_1,s_2) &= \zeta_{10}s_1+\zeta_{01}s_2 + \zeta_{00},
	\end{aligned}
\end{equation}
where \(\zeta_{ij} \in \R\), and $\psi_1$, $\psi_0$ $\psi_{\frac{1}{2}}$ are linear functions in $s=(s_1,s_2)\in\R^2$. Justifications are in~\cite{KUMAR-19} and details in~\cite{KARMARKAR-24}.

To complete Step 3 and obtain the conjugate of the entire PLQ function, we work in two stages: Step 3a given $s_{i,j,k}$, compute $s_{i,k}$, and Step 3b given $s_{i,k}$, compute $s_{k}$.

\subsubsection{Step 3a}
We divide the task into three steps: (i) compute intersection of the domains, (ii) compute maximum and if needed split the domain, and (iii) merge adjacent regions if they have the same conjugate expressions.

\paragraph{Step 3a(i)}
Let us first consider the domains of the conjugates corresponding to the first two convex functions in the $i^{th}$ piece. The domain of the first conjugate is $\cup_{k_1=1 \dots n_{s_{i,1}}} \dom  s_{i,1,k_1} $ while $\cup_{k_2=1 \dots n_{s_{i,2}}} \dom   s_{i,2,k_2}$ is the domain of the second conjugate of the $i^{th}$ piece. We need to find the intersection $ \dom s^t_{i,k} = \dom s_{i,1,k_1} \cap \dom s_{i,2,k_2}$, for all $k_1$, $k_2$.

\begin{figure}
	\centering
	\subfigure[$\dom s_{i,1,k_1}$, $k_1 = 1\dots 6$.]{ \label{fig:s11}
		\centering
		\includegraphics[height=127px]{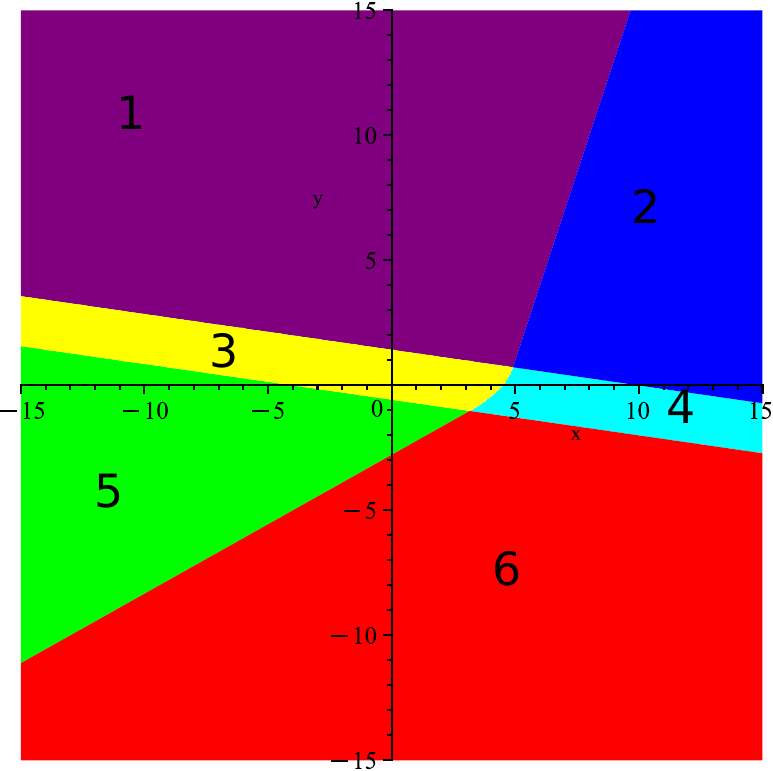}}
	\subfigure[$\dom s_{i,2,k_2}, k_2 = 1 \dots 3.$]{\label{fig:s12}
		\centering
		\includegraphics[height=130px]{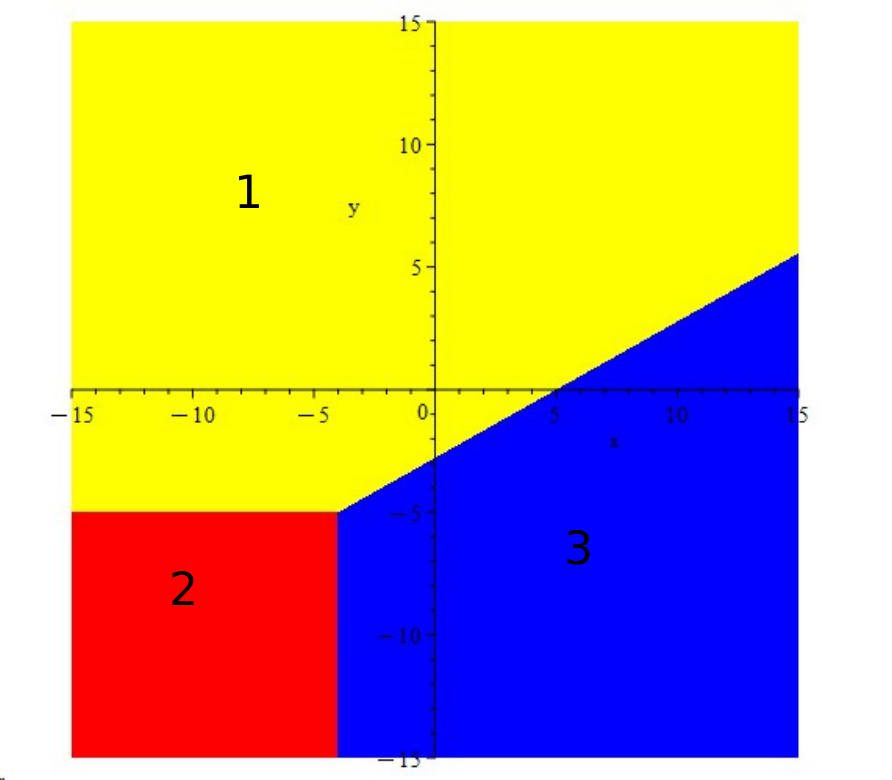}}
	\\
	\subfigure[$\dom s^t_{i,k}=\dom s_{i,1,k_1} \cap \dom s_{i,2,k_2}.$]{
		\centering
		\includegraphics[height=150px]{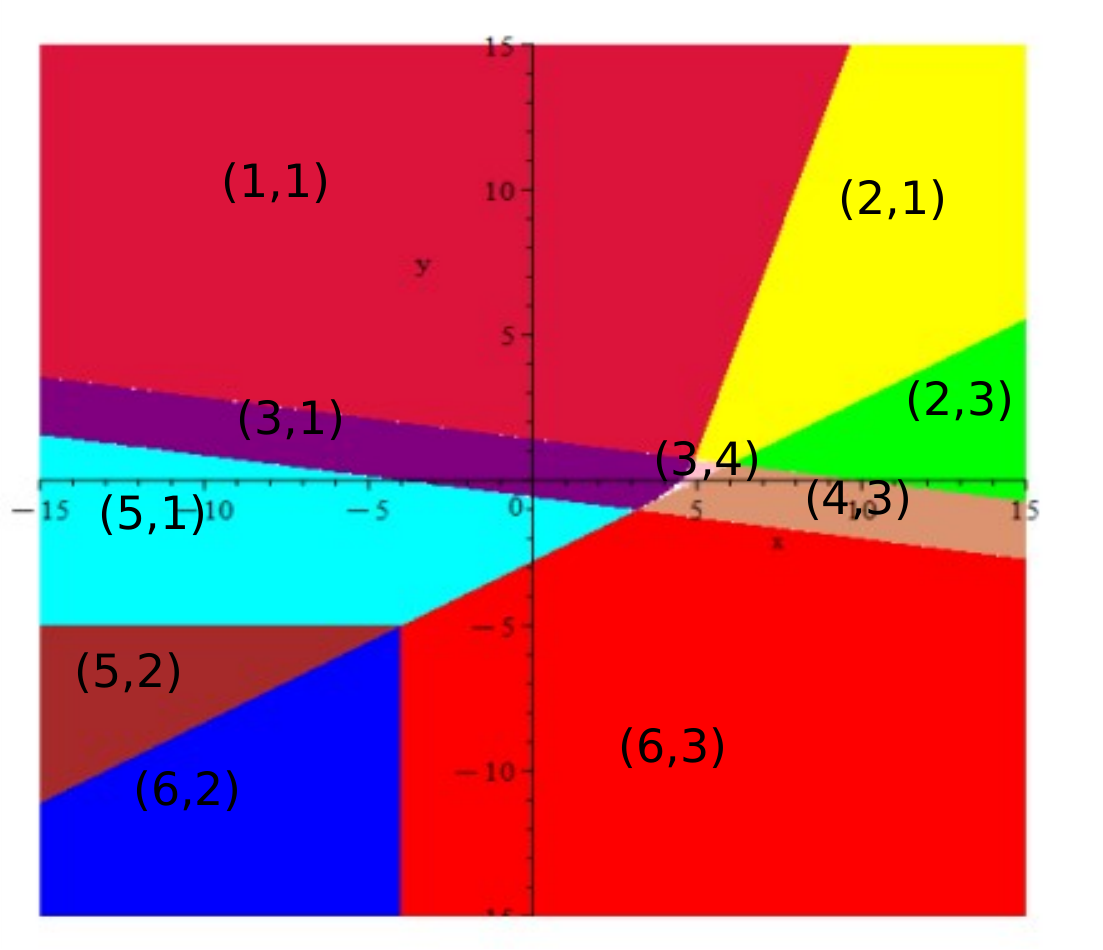}
		\label{fig:s1t}
	}
	\caption{Intersection of domains of two conjugates.}
	\label{fig:s3a1}
\end{figure}

A simple algorithm is to exhaustively generate all pairs between the two sets of regions and check if the intersection of the interior of the two regions is nonempty. When we find a nonempty intersection, we store it in the output list along with the two functions which were defined on the two regions we were intersecting. So we are basically storing nonempty intersections along with two functions, as $(\dom s^t_{i,k},[s_{i,k,1},s_{i,k,2}]) = (\dom s_{i,1,k_1} \cap \dom s_{i,2,k_2},[s_{i,1,k_1},s_{i,2,k_2}])$. The output $\dom s^t_{i,k}$ is a list of nonempty domains covering the entire $\R^2$ plane as shown in Figure~\ref{fig:s3a1}. 

The above algorithm considers all pairs so it has a complexity of $O(n^2)$, where $n=\max(n_{s_{i,1}},n_{s_{i,2}})$, but it is possible to reduce it to $O(n)$. When we consider regions that are neighbors, we get nonempty intersections only when the original convex envelopes have a common vertex or edge. So a propagation algorithm similar to~\cite{HAQUE-18} provides a linear number of nonempty intersections. Computing the intersection of two convex polytopes is linear in the number of vertexes~\cite{CHAN-16} resulting in a linear-time algorithm. Practically, starting with a triangulation, we have a bounded number of vertices. Computing a triangulation for a convex polyhedral set takes linear time since triangulating such sets only require joining non-neighbor vertexes from any arbitrary vertex~\cite[p. 49]{BERG-08}. Hence, our approach runs in linear-time.

At this stage we are finding the intersection of domains of the original conjugates. These domains are parabolic subdivisions~\cite{KUMAR-19}. After computing the intersection we again get a parabolic subdivision.

\begin{proposition}\label{prop:1}
	The output of Step 3a(i) is a parabolic subdivision.
\end{proposition}

\begin{proof}
	Each $\dom s_{i,j,k}$ is a parabolic subdivision, which is the nonempty interior of a set of parabolic inequalities. So each domain is of the form $\dom s_{i,j,k} = \{ (x,y),a_lx^2+b_lxy+c_ly^2+d_lx+e_ly+f_l < 0, l=1\dots K\}$ where $b_l^2-4a_lc_l = 0$, and $K \in \N $ is the number of inequalities defining $\dom s_{i,j,k}$. When we find the intersection of two such domains, we only store the regions with nonempty interiors. These are again defined by parabolic inequalities of the same form, hence the subdivision is still parabolic.
\end{proof}

\begin{figure}
	\centering
	\subfigure[$\co (q_i + I_{\dom q_i})= \min_{j=1 \dots 2}(r_{ij} + I_{\dom r_{ij}})$. ]{
		\centering
		\includegraphics[width=0.45\textwidth]{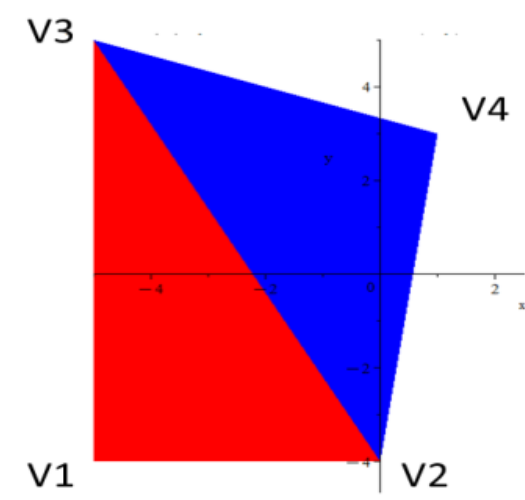}\label{f:intersection_Ona}}
	\\
	\subfigure[$\dom (s_{i,1,k}).$ ]{
		\centering
		\includegraphics[width=0.35\textwidth]{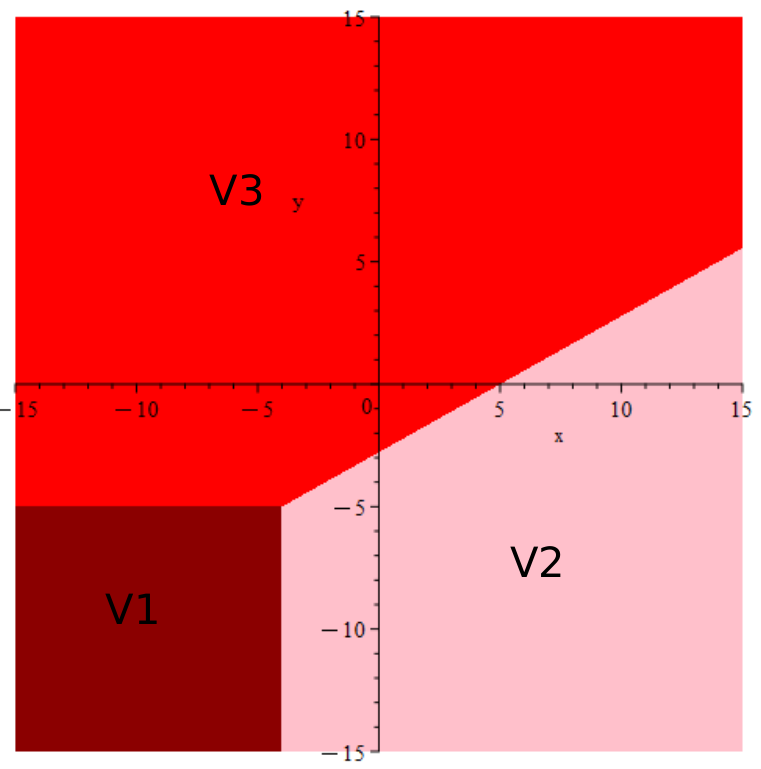}\label{f:intersection_Onb}}
	\subfigure[$\dom (s_{i,2,k}).$ ]{
		\centering
		\includegraphics[width=0.35\textwidth]{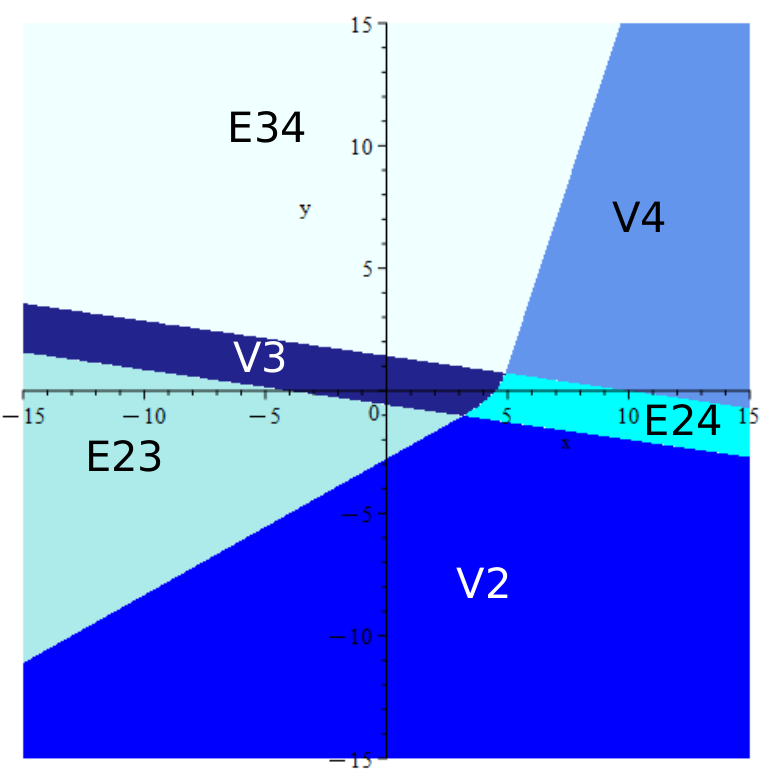}\label{f:intersection_Onc}}
	\\
	\subfigure[$\dom (s_{i,1,k})\cap\dom (s_{i,2,k}).$ ]{
		\centering
		{\includegraphics[width=0.45\textwidth]{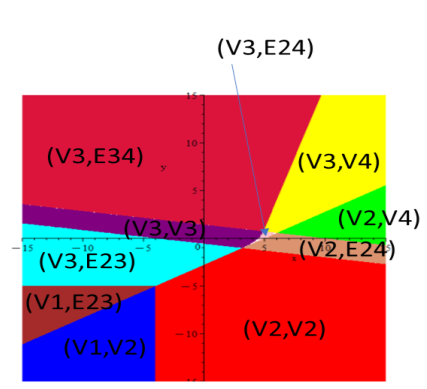}\label{f:intersection_Ond}}}
	\caption{Computing $\dom (s_{i,1,k})\cap\dom (s_{i,2,k})$ using $\dom (r_{i,1})\cap\dom (r_{i,2}).$ }
	\label{fig:intersection_On}
\end{figure}

We illustrate Step 3a with an example where $V$ denotes a vertex and $E$ denotes the relative interior of an edge. Figure~\ref{f:intersection_Ona}, shows the convex polyhedral subdivision of the domain of one piece as obtained in Step 1. As a result of Step 1 we get $r_{i1} = \co (V_1,V_2,V_3)$ and $r_{i2} = \co (V_2,V_3,V_4)$. 
Figure~\ref{f:intersection_Onb} and Figure~\ref{f:intersection_Onc} show  $\dom(s_{i1})$ and $\dom(s_{i2})$ respectively. From Figure~\ref{f:intersection_Ona} we observe that $\dom(r_{i1}) \cap \dom(r_{i2}) = \{V_2, V_3, E_{23}\}$. The set $\dom(s_{i1})$ has pieces corresponding to $R_1 = \{V_1,V_2, V_3\}$,  whereas $\dom(s_{i2})$  has pieces corresponding to $R_2= \{V_2,V_3, V_4, E_{23}, E_{24}, E_{34}\}$.
In Figure~\ref{f:intersection_Ond}, the nonempty intersections always have an entity corresponding to $R_1 \cap R_2$.

\paragraph{Step 3a(ii)}
From Step 3a(i) we get a list $\dom s^t_{i,k}$, along with two functions $s_{i,k,1}$ and $s_{i,k,2}$ defined over each $\dom s^t_{i,k}$. We now find the maximum function $s^t_{i,k} = \max\{s_{i,k,1},s_{i,k,2}\}$ over the given $\dom s^t_{i,k}$. In doing so it is possible that $\dom s^t_{i,k}$ is further divided into two regions one where $s_{i,k,1}$ is the maximum, and the other where the maximum is $s_{i,k,2}$; or we get the same region with either $s_{i,k,1}$ or $s_{i,k,2}$ as the maximum. 

\begin{figure}
    \centering
	\includegraphics[width=12cm]{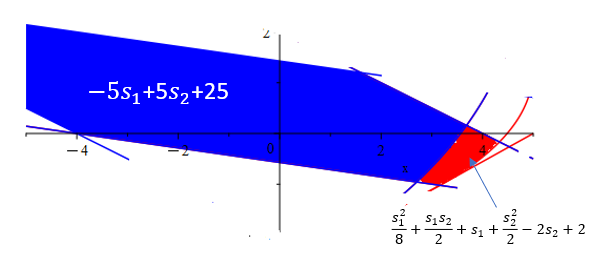}
	\caption{Further division by a parabolic curve as explained in Example \ref{division}.}
	\label{fig:maxP}
\end{figure}

\begin{example} \label{division}
	The blue region in Figure~\ref{fig:maxP} is denoted by $B$ while $R$ is the red region; they are defined as 
	\begin{multline*}
			B=\{(u,v):-u -7v -4 \le 0 ,
			u +7v -10 \le 0, 
			-u -2v -4 \le 0 ,\\
			u +2v -4 \le 0, 
			48u - 56v +4uv +u^2 +4v^2 - 184 \le 0\};
	\end{multline*}
	\begin{multline*}
			R = \{(u,v):-u -7v -4 \le 0 ,
			148u - 196v +(u +7v)^2 -684 \le 0 ,
			u +2v -4 \le 0,\\
			56v -48u - 4uv -u^2 -4v^2 +184 \le 0 \}.
	\end{multline*}
	
	Functions $f_1(u,v)=-5u+5v+25$ and $f_2(u,v)=u^2/8+uv/2+u+v^2/2-2v+2$ are defined on $B \cup R$. Computing $\max\{f_1,f_2\}$ over $B \cup R$ gives a further subdivision along $f_1-f_2=0$; see Figure~\ref{fig:maxP}. Function $f_1$ is the maximum over $B$ while $f_2$ is over $R$.
\end{example}

At this point, when we find the maximum of the two functions defined on the intersecting domain, we get a further subdivision with an inequality  $f_1-f_2\le 0$ for one region and $f_2-f_1\le 0$ for the other. Thus in the maximum we can get a subdivision which is the difference of the two conjugate functions, $f_1$ and $f_2$. From \cite[Theorem 4.24, 4.27]{KUMAR-19} we get the conjugate functions of the forms~\eqref{eq:c2}. In order to find the possible subdivisions, we enumerate all the pairs, of the conjugate functions, $(f_1,f_2)$ that give us a further subdivision of the domain. We first give propositions for each possible case and the partition it can give us and then put everything together in Theorem~\ref{thm:subdiv}.

\begin{proposition} \label{div_linlin}
	Let the two functions defined on the parabolic domain be linear. If we find the $\max\{f_1,f_2\}$, we can get a further subdivision of the domain. This division is given by a linear inequality hence the subdivision remains parabolic.
\end{proposition}

\begin{proof}
	The domain is split by the equality, $f_1-f_2=0$ which is linear. 
\end{proof}

In order to get the subdivision when one conjugate expression is quadratic, we first need to show that the quadratic expression $g_q$ obtained in~\eqref{eq:c2} is always parabolic.

\begin{proposition}\label{thm_pexpr}
	The quadratic expression $g_q$ obtained in~\eqref{eq:c2} is parabolic.        
\end{proposition}

\begin{proof}
	The conjugate expression obtained is a quadratic only in a certain condition. In Step 1, when we find the convex envelope of each piece, we can obtain a bivariate rational function, $\frac{\xi_1^2(x)}{\xi_2(x)}+\xi_0(x)$, where $\xi_0(x)=\xi_{01}x_1+\xi_{02}x_2+\xi_{00}$, $\xi_1(x)=\xi_{11}x_1+\xi_{12}x_2+\xi_{10}$ and $\xi_2(x)=\xi_{21}x_1+\xi_{22}x_2+\xi_{20}$. Then in Step 2, we find the conjugate expression of this rational function. Using the method given in~\cite{KUMAR-19}, when we find the expression on an edge ($y=m x+q$) of the polytope and $\xi_{21}+m\xi_{22}=0$, we get a quadratic expression. The quadratic expression obtained is 
	\begin{equation}\label{eq:pq}
		f^*(s)=\zeta_{11}s_1^2+\zeta_{12}s_1s_2+\zeta_{22}s_2^2+\zeta_{10}s_1+\zeta_{01}s_2+\zeta_{00},    
	\end{equation}
	where
	\[\zeta_{11}=-\frac{(\xi_{11}\gamma_{10}+m\xi_{12}\gamma_{10})^2}{\xi_{20}+q\xi_{22}}+\gamma_{10},\]
	\[\zeta_{12}=-\frac{2(\xi_{11}\gamma_{01}+m\xi_{12}\gamma_{01})(\xi_{11}\gamma_{10}+m\xi_{12}\gamma_{10})}{\xi_{20}+q\xi_{22}}+\gamma_{01}+m\gamma_{10},\]
	\[\zeta_{22}=-\frac{(\xi_{11}\gamma_{10}+m\xi_{12}\gamma_{10})^2}{\xi_{20}+q\xi_{22}}+\gamma_{10}m,\]
	\begin{equation}
		\begin{split}
			\zeta_{10}=-\frac{2(\xi_{11}\gamma_{10}+m\xi_{12}\gamma_{10})(\xi_{10}+\xi_{11}\gamma_{00}+\xi_{12}(q+m\xi_{00}))}{\xi_{20}+q\xi_{22}}+\\
			\gamma_{00}-m\xi_{02}\gamma_{10}-\xi_{01}\gamma_{10},  
		\end{split}
	\end{equation}

\begin{equation}
	\begin{split}
		\zeta_{01}=-\frac{2(\xi_{11}\gamma_{01}+m\xi_{12}\gamma_{01})(\xi_{10}+\xi_{11}\gamma_{00}+\xi_{12}(q+m\xi_{00}))}{\xi_{20}+q\xi_{22}}+\\m\gamma_{00}-m\xi_{02}\gamma_{01}-\xi_{01}\gamma_{01}+q,\end{split}
\end{equation}
and
\[\zeta_{00}=-\frac{(\xi_{10}+\xi_{11}\gamma_{00}+\xi_{12}(q+m\xi_{00}))^2}{\xi_{20}+q\xi_{22}}-\xi_{00}-\xi_{01}\gamma_{00}-\xi_{02}(m\gamma_{00}+q).\]
Now substituting the above values, we obtain
\[\zeta_{12}^2 - 4\zeta_{11}\zeta_{22}=0.\]
Hence~\eqref{eq:pq} is parabolic.
\end{proof}

The MATLAB code to verify this proof is given in Section~\ref{a:proof_them_pexpr}.

\begin{lemma} \label{div_linquad}
If $f_1$ is a parabolic quadratic function and $f_2$ is linear, computing $\max\{f_1,f_2\}$ may result in a division by the curve $f_1-f_2=0$. In that case, the curve $f_1-f_2=0$ is parabolic.
\end{lemma}

\begin{proof}
Let the parabolic expression defined be $f^1(s)=\zeta_{11}^1s_1^2+\zeta_{12}^1s_1s_2+\zeta_{22}^1s_2^2+\zeta_{10}^1s_1+\zeta_{01}^1s_2+\zeta_{00}^1 $ and the linear function be $f^2(s)=\zeta_{10}^2s_1+\zeta_{01}^2s_2+\zeta_{00}^2$.
The new subdivision is given by $f^1(s)-f^2(s) = \zeta_{11}^1s_1^2+\zeta_{12}^1s_1s_2+\zeta_{22}^1s_2^2+(\zeta_{10}^1-\zeta_{10}^2)s_1+(\zeta_{01}^1-\zeta_{01}^2)s_2+(\zeta_{00}^1-\zeta_{00}^2)=0$.
As $f^1(s)$ is parabolic, we have $\zeta_{12}^2 - 4\zeta_{11}\zeta_{22}=0$, hence $f^1(s)-f^2(s)$ is also parabolic.
\end{proof}
This is shown in Figure~\ref{fig:maxP}.

\begin{theorem}\label{thm:subdiv}
The maximum of two conjugates is a piecewise function that admits a parabolic subdivision.
\end{theorem}

\begin{table}[tbph]
\centering
\caption{Possible divisions of domain obtained in Step 3a (ii).
	\label{table:observedivision}}
\begin{tabular}{@{}llr@{}} \toprule 
	$s_{i,1,k}$ & $s_{i,2,k}$ & $s_{i,1,k}-s_{i,2,k}$ \\\toprule 
	Linear & Linear & Linear \\ 
	Linear & Parabolic & Parabolic \\ 
	\bottomrule 
\end{tabular}
\end{table}

\begin{proof}
When we find the convex envelope of each piece in Step 1, we do not get rational functions in adjacent domains. They only arise when the region is divided into three regions. The middle region has the rational function, while the adjacent regions have linear expressions. 
In Step 2, only the rational functions restricted to an edge give rise to quadratic expressions.
In Step 3, when we compare two functions we always get one of them as linear. Nonempty intersections are obtained only when the primal has a common edge. Hence the two functions we compare always have one linear function. Thus the additional subdivision is one of the forms given in Table~\ref{table:observedivision}. The proofs that these subdivisions are parabolic are given in propositions \ref{div_linlin}, and \ref{div_linquad}.
\end{proof}

\begin{figure}
	\centering
	\includegraphics[width=12cm]{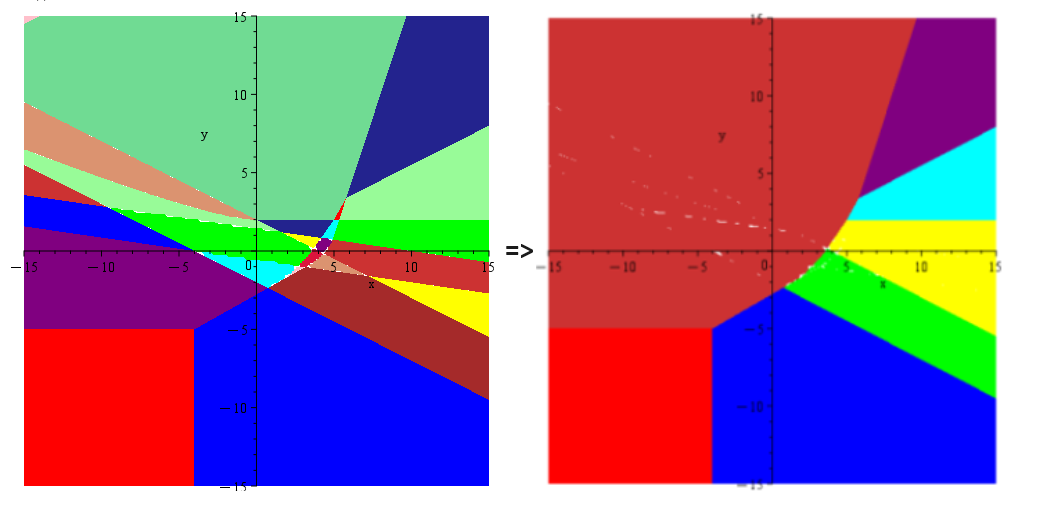}
	\caption{Merging adjacent pieces with the same function.}
	\label{fig:merge}
\end{figure}

\paragraph{Step 3a(iii)}
At this point we  computed the maxima and obtained the conjugate as a piecewise function. This can be further simplified by merging adjacent regions where functions are the same as shown in Figure~\ref{fig:merge}. Coefficients of the functions involved are stored as rational numbers, so we can determine two identical functions without floating point considerations. 

For every piece of the PLQ function, we have multiple conjugates, $r_{ijk}$. So far, we have given a method to compute the maximum of two conjugates. We can now add one conjugate at a time to get the conjugate of the entire piece. We first copy the first conjugate into the output. We then iterate over the remaining pieces, performing steps 3a(i) and 3a(ii) and finally get the conjugate of the entire piece.

\paragraph{Division}

\begin{proposition}\label{p:maxParabolic}
    The function $\max_j(r_{ijk})$ admits a parabolic subdivision. 
\end{proposition}
	
\begin{proof}
    By Proposition~\ref{prop:1} and Theorem~\ref{thm:subdiv}, the maximum of two conjugates is a parabolic subdivision. In computing the convex envelope of a piece only one division can have a rational function defined on it. In computing the conjugate, only the rational function would give us a quadratic form. Thus the maximum of multiple divisions is a parabolic subdivision.
\end{proof}

\begin{figure}
	\centering
	\subfigure[$s_{1,1,k}.$ ]{
		\centering
		{\includegraphics[width=.45\textwidth]{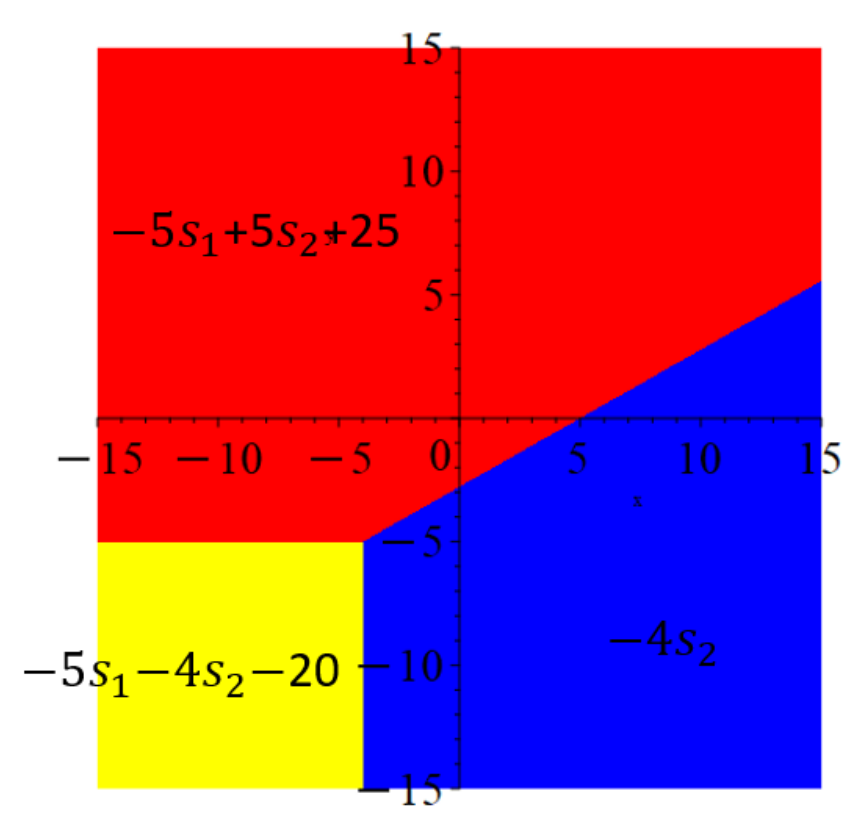}}}
	\subfigure[$s_{1,2,k}.$ ]{
		\centering
		{\includegraphics[width=.418\textwidth]{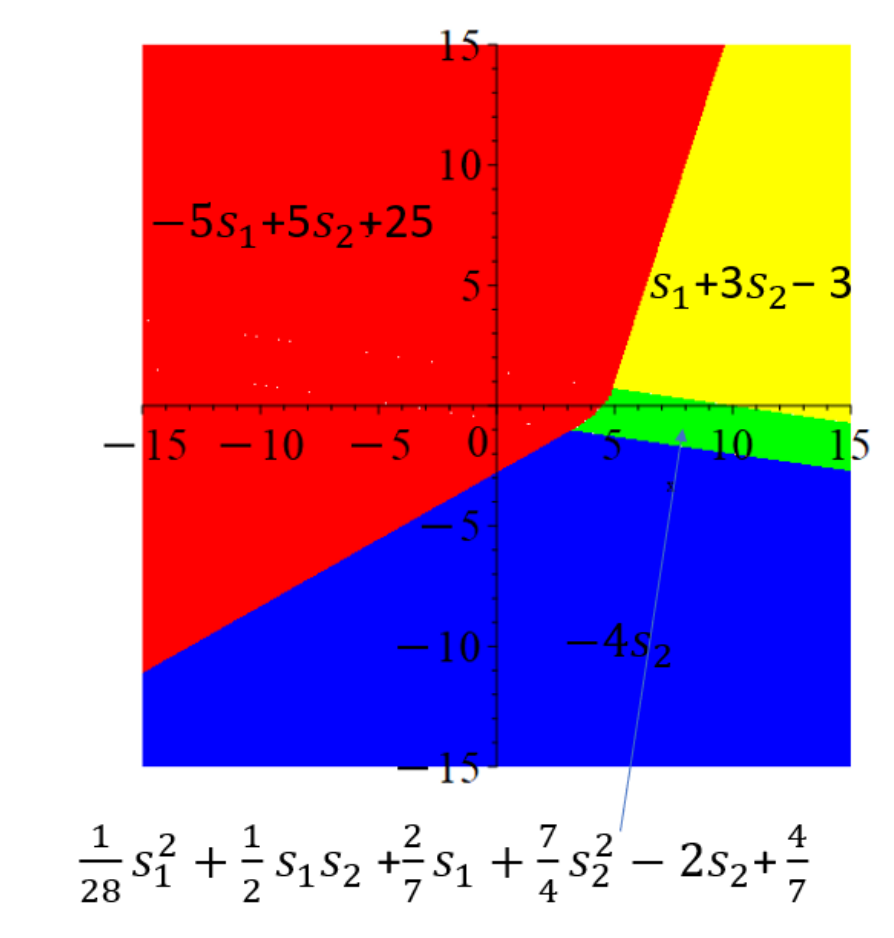}}}
	\\
	\subfigure[$\max\{ s_{1,1,k}s_{1,2,k}\}.$ ]{
		\centering
		\includegraphics[width=.7\textwidth]{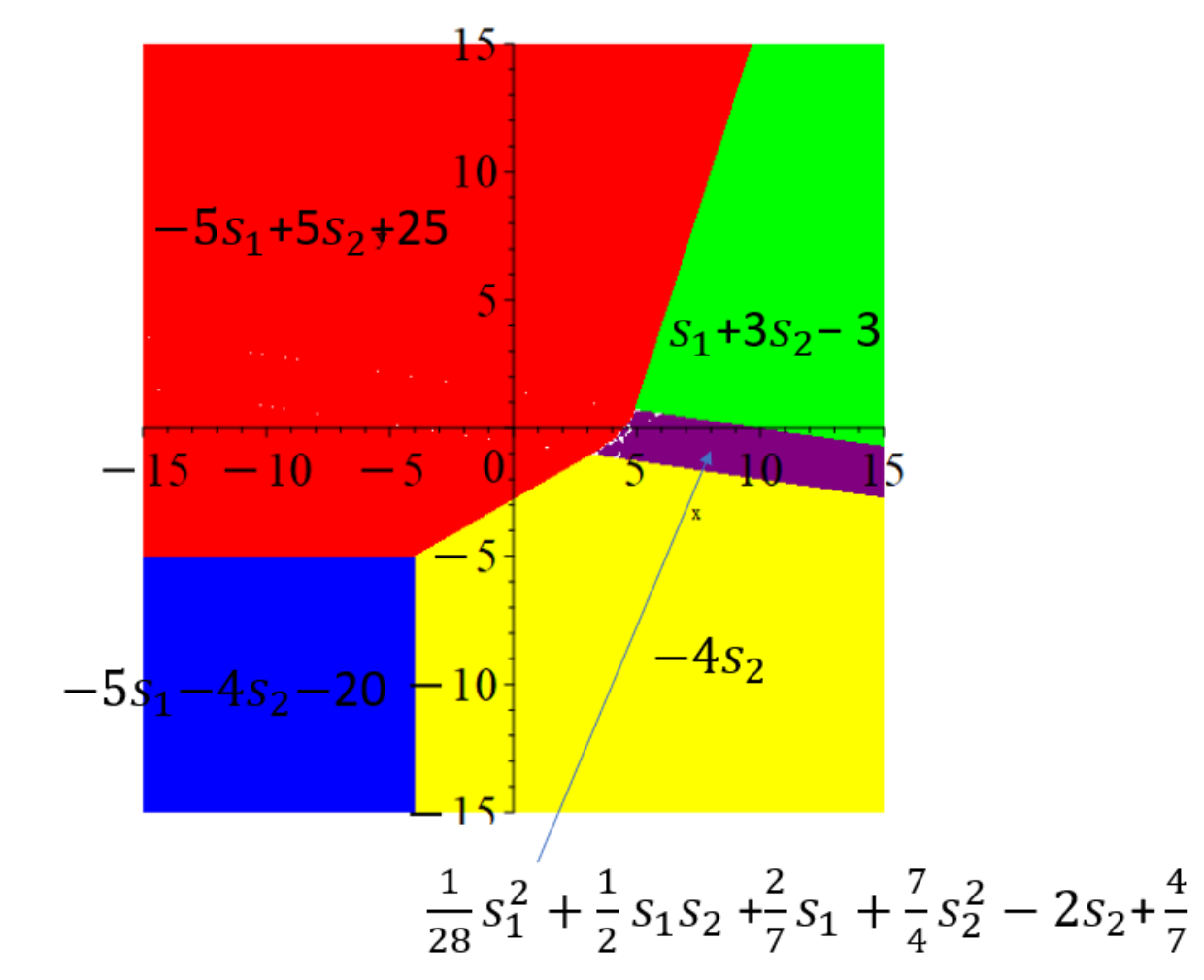}}
	\caption{Maximum Conjugate for piece 1.}
	\label{fig:conjugate1}
\end{figure}

\begin{example}\label{example1} 
	In the example in Figure~\ref{fig:example1-all}(a), we had two pieces. On finding the convex envelope in Step 1 we obtained convex functions defined over a polyhedral subdivision, refer to Figure~\ref{fig:example1-all}(b-c). Corresponding to these, in Step 2 we got two piecewise conjugates as shown in Figure~\ref{fig:example1-all}(d-g). The maximum of the conjugates for Piece 1, is shown in Figure~\ref{fig:conjugate1}.
\end{example}

\begin{figure}
	\centering
	\subfigure[$s_{2,1,k}.$ ]{
		\centering
		\includegraphics[width=.45\textwidth]{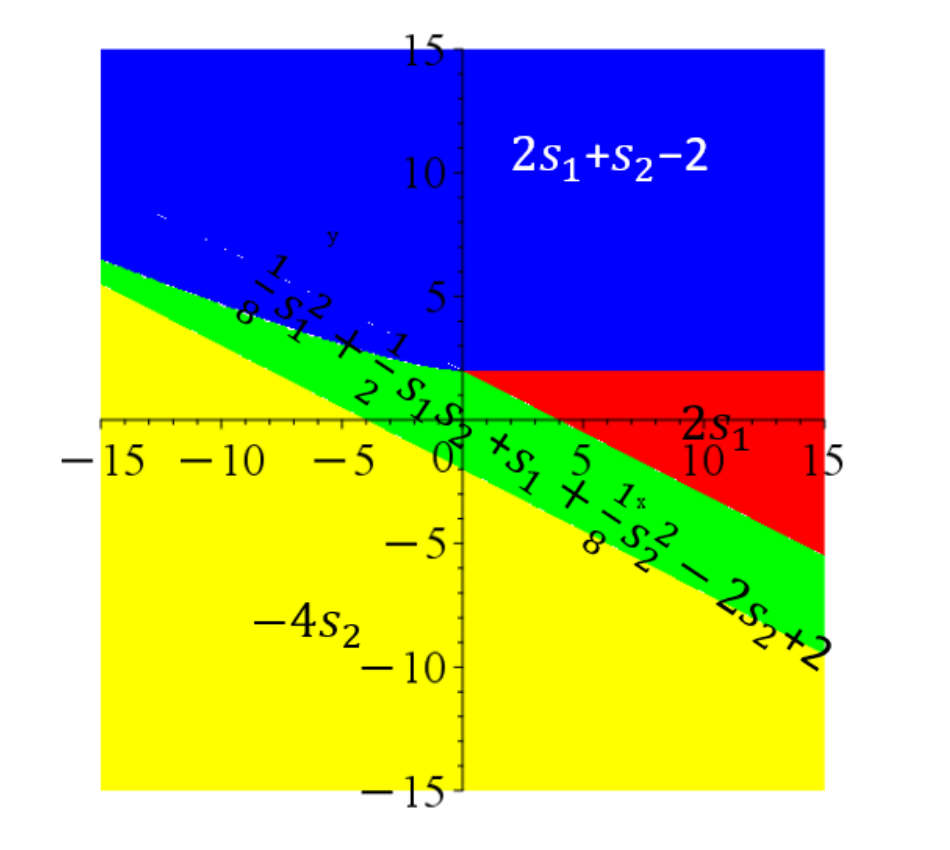}}
	\subfigure[$s_{2,2,k}.$ ]{
		\centering
		\includegraphics[width=.45\textwidth]{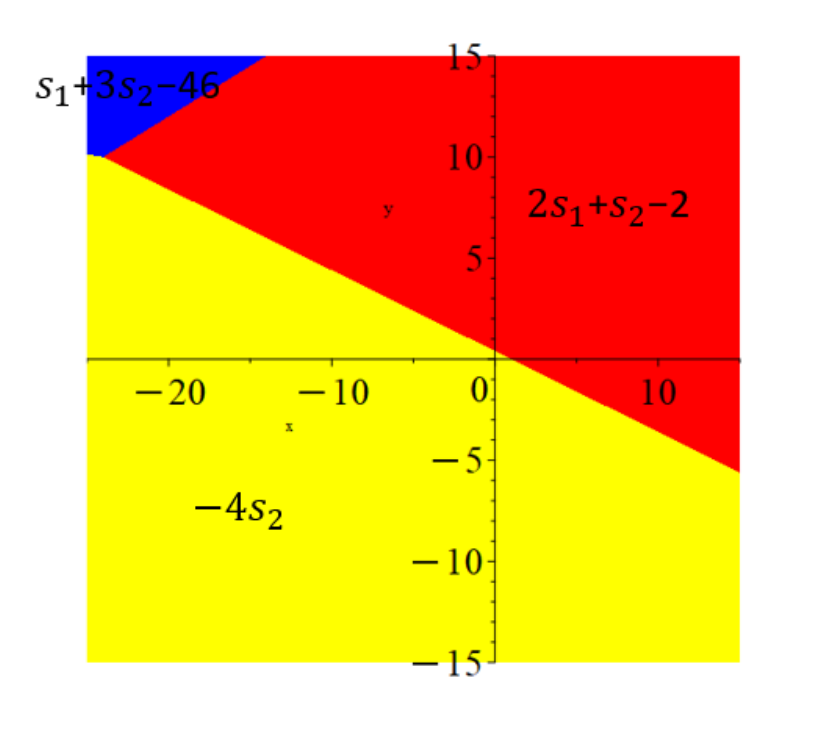}}
	\\
	\subfigure[$\max\{ s_{2,1,k}s_{2,2,k}\}.$ ]{
		\centering
		\includegraphics[width=.6\textwidth]{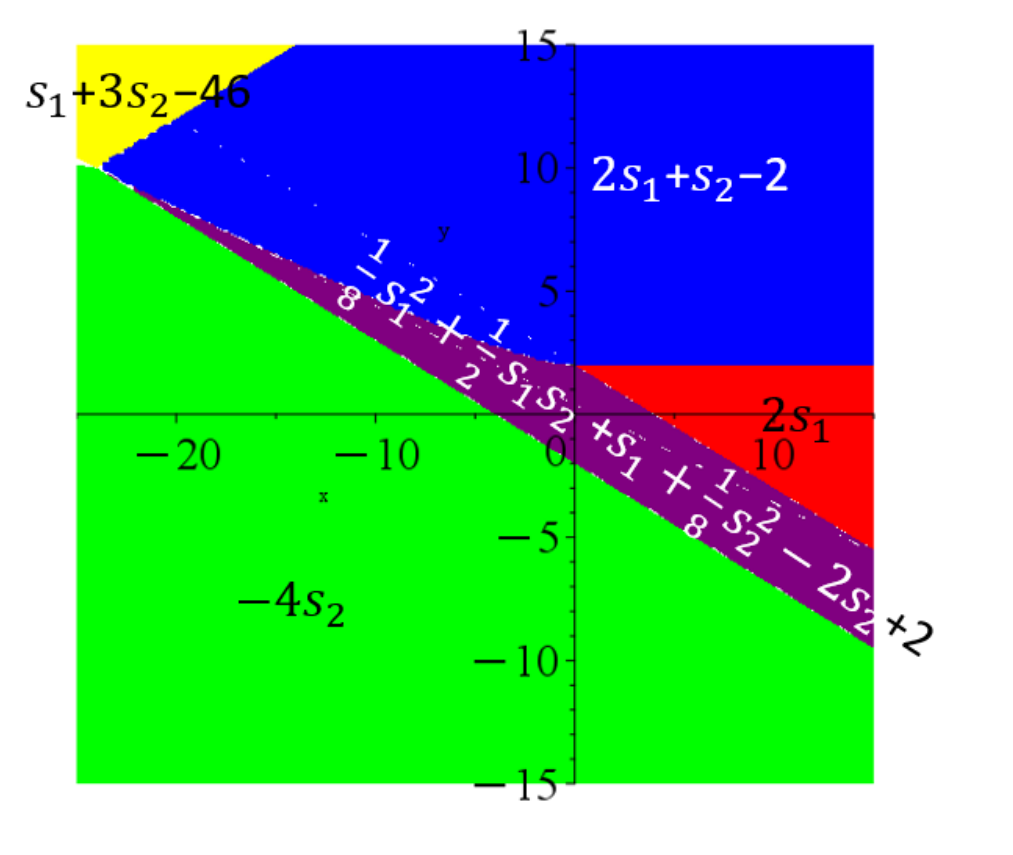}}
	\caption{Maximum Conjugate for piece 2.}
	\label{fig:conjugate2}
\end{figure}

The same process is repeated for Piece 2 shown in Figure~\ref{fig:1}, and the maximum conjugate for this piece is shown in Figure~\ref{fig:conjugate2}.

\subsubsection{Step 3b}
In Step 3a, corresponding to each piece of the PLQ function, we obtained a maximum conjugate which is the maximum of all the conjugates for that piece. 
In this stage we find the maximum over all these maximums to get one maximum for the entire PLQ function. The steps followed to compute this are very similar to those followed in Step 3a.   
In Step 3b(i) we find the intersection of the domains, $ \dom s^t_{k} = \dom s_{1,k} \cap \dom s_{2,k} $  and in Step 3b(ii), we find the maximum of two functions over the same domain as $s^t_{k}=\max\{s_{k,1},s_{k,2}\}$ over $\dom s^t_{k}$ and then merge to get the maximum over the first two pieces. We then iterate Step 3b(i) and Step 3b(ii), over the remaining pieces, $\dom s^t_{k}=\dom s^t_{k} \cap \dom s_{j,k}$, $s^t_{k}=\max\{s^t_{k},s_{j,k}\}$ and merge to finally get the maximum conjugate \(s_{k} = s^t_{k}\), which is the conjugate of the PLQ function.

\paragraph{Divisions}
In this step we are finding the maximum between two pieces. In Step 3b(i), we find the intersection of parabolic subdivisions and hence get parabo\-lic subdivisions as given in Proposition~\ref{prop:1}. In Step 3b(ii), we are finding the maximum of two functions over a parabolic domain.
The pieces are adjacent to each other in the primal and Proposition~\ref{p:maxParabolic} shows that we would get a parabolic subdivision as if we had considered both pieces together as one piece.

\begin{figure}
	\centering
	\subfigure[$s_{1,k}.$ ]{
		\centering
		\includegraphics[width=0.45\textwidth]{example1_conjmax1-eps-converted-to.pdf}}
	\subfigure[$s_{2,k}.$ ]{
		\centering
		\includegraphics[width=0.45\textwidth]{example1_conjmax2-eps-converted-to.pdf}}\\
	\subfigure[$\max\{ s_{1,k}s_{2,k}\}.$ ]{
		\centering
        \includegraphics[width=0.65\textwidth]{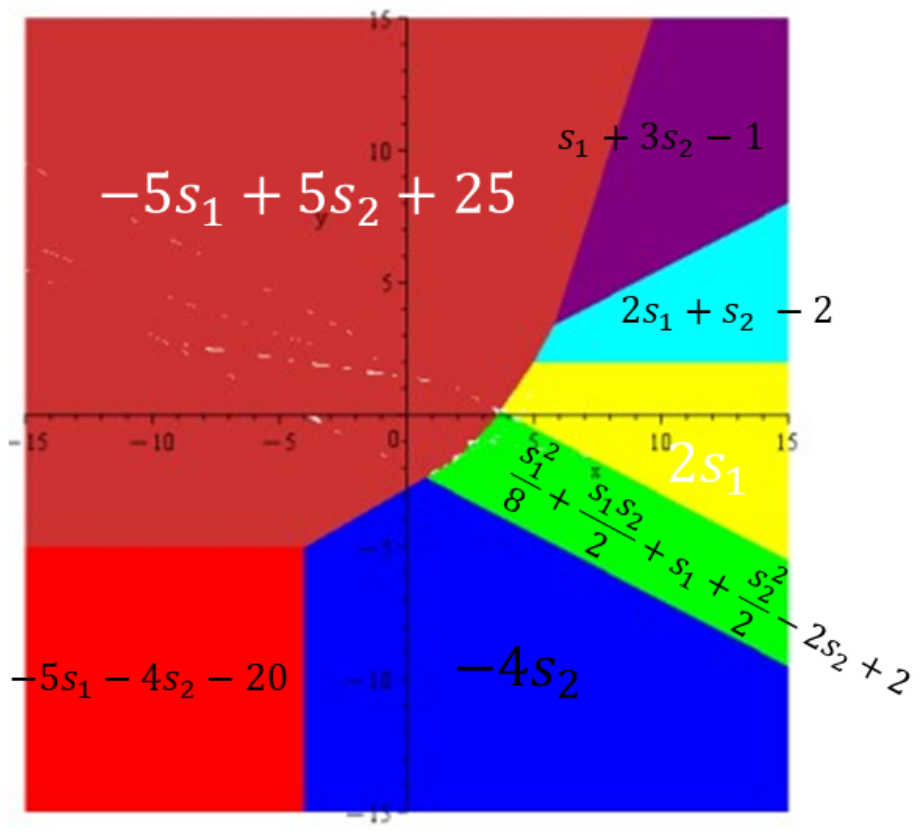}}
	\caption{Maximum Conjugate over all pieces.}
	\label{fig:conjugate}
\end{figure}

\begin{example}
	Again, consider the example given in Figure~\ref{fig:1}. We had two pieces. Corresponding to each piece, in Step 3a we got the maximum conjugate of each piece which were shown in figures~\ref{fig:conjugate1} and~\ref{fig:conjugate2}. We now find the maximum over these maximums to get the conjugate of the PLQ function as shown in Figure~\ref{fig:conjugate}. 
\end{example}

\begin{figure}
	\centering
	\subfigure[ $q_{1}.$ ]{
		\centering
		\includegraphics[height=80px]{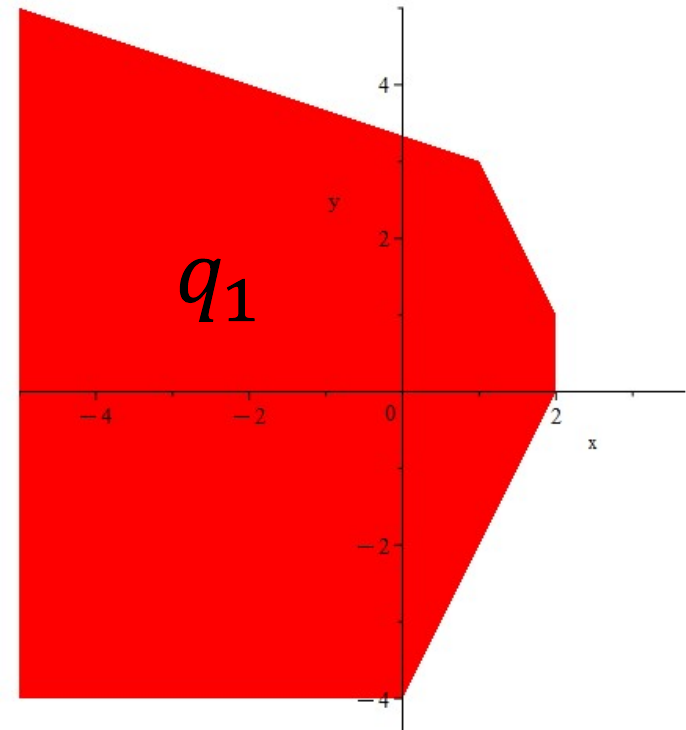}}	\\
	\subfigure[$r_{1,j}.$ ]{
		\centering		
        \includegraphics[height=120px]{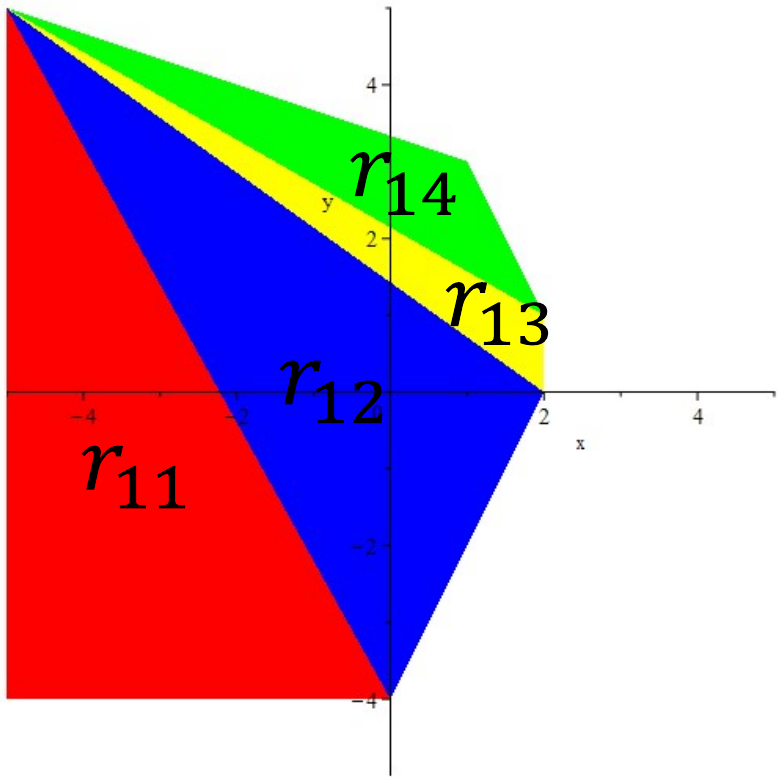}}	\\
	\subfigure[$s_{1,1,k}.$ ]{
		\centering
        \includegraphics[height=80px]{latexpics_conjugate1_1T.pdf}}
	\subfigure[$s_{1,2,k}.$ ]{
		\centering
		\includegraphics[height=80px]{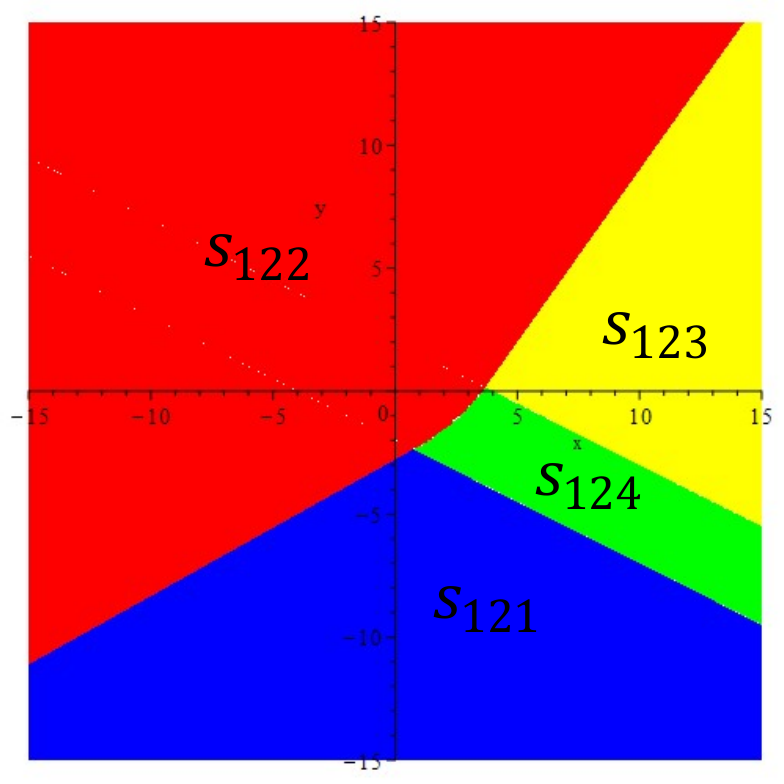}}
	\subfigure[$s_{1,3,k}.$ ]{
		\centering
        \includegraphics[height=80px]{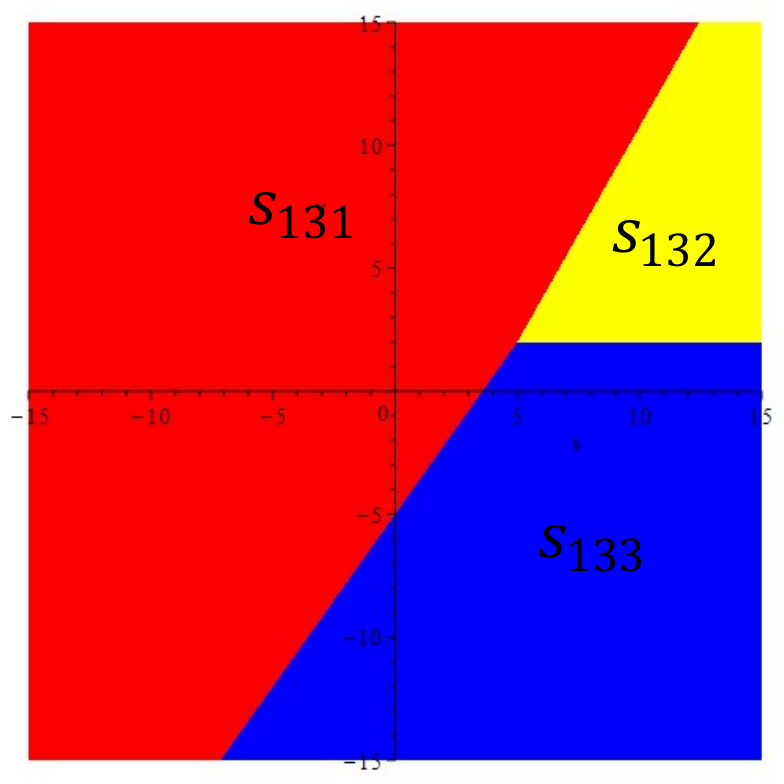}}
	\subfigure[$s_{1,4,k}.$ ]{
		\centering
        \includegraphics[height=80px]{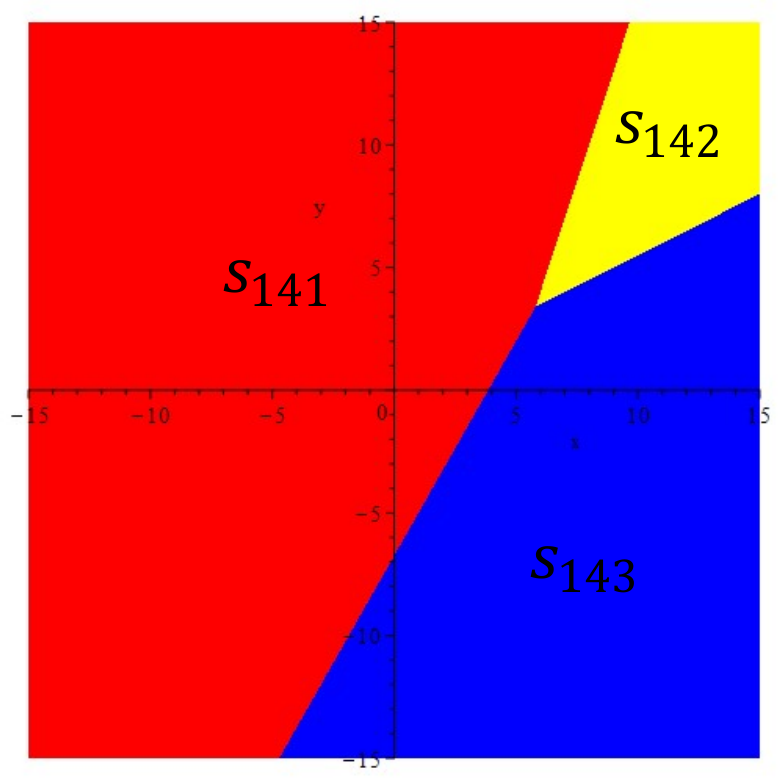}}	\\
	\subfigure[$s_{k}.$ ]{
		\centering
        \includegraphics[height=165px]{latexpics_maxT.pdf}}
	
	\caption{Conjugate of Example \ref{ex:overall2}.(The functions are enumerated in the Table \ref{table:allsteps2})}
	\label{fig:conjugateMax2}
\end{figure}

\begin{example} \label{ex:overall2}
	Here we refer to the PLQ function shown in Figure~\ref{fig:1} but instead of dividing it into two pieces we take just one piece as illustrated in Figure~\ref{fig:conjugateMax2} and use this to test that the final output for both examples is the same. In this example as there is only one piece, Step 3b is not required.
\end{example}

\begin{table}[tbph]
\centering
\caption{Function definitions for Figure \ref{fig:conjugateMax2}.}
	\label{table:allsteps2}
\begin{tabular}{@{}llr@{}} \toprule 
    Function name & Expression\\
    \midrule
	$q_1$ &$xy$ \\
    $r_{11}$ &$-4x-5y-20$  \\ 
    $r_{12}$ &$(30x - 5y + 4xy+10x^2+5y^2-100)/(2*x - y + 15)$  \\ 
    $r_{13}$ &$5x + 2y - 10$  \\ 
    $r_{14}$ &$\frac{29x}{5} + \frac{17y}{5}-13$  \\ 
    $s_{111}$ &$-5s_1+4s_2-20$  \\ 
    $s_{112}$ &$-5s_1+5s_2+25$  \\ 
    $s_{113}$ &$-4s_2$  \\ 
    $s_{121}$ &$-4s_2$  \\ 
    $s_{122}$ &$-5s_1+5s_2+25$  \\ 
    $s_{123}$ &$2s_1$  \\ 
    $s_{124}$ &$\frac{1}{8}s_1^2+\frac{1}{2}s_1s_2+s_1+\frac{1}{8}s_2^2-2s_2+2$  \\ 
    $s_{131}$ &$-5s_1+5s_2+25$  \\ 
    $s_{132}$ &$s_1+3s_2-3$  \\ 
    $s_{133}$ &$2s_1$  \\ 
    $s_{141}$ &$-5s_1+5s_2+25$  \\ 
    $s_{142}$ &$s_1+3s_2-3$  \\ 
    $s_{143}$ &$2s_1+s_2-2$  \\ 
    $s_{1}$ &$-5s_1+4s_2-20$  \\ 
    $s_{2}$ &$-5s_1+5s_2+25$  \\ 
    $s_{3}$ &$s_1+3s_2-3$  \\ 
    $s_{4}$ &$2s_1+s_2-2$  \\ 
    $s_{5}$ &$2s_1$  \\ 
    $s_{6}$ &$\frac{1}{8}s_1^2+\frac{1}{2}s_1s_2+s_1+\frac{1}{8}1s_2^2-2s_2+2$  \\ 
    $s_{7}$ &$-4s_2$  \\ 
	\bottomrule 
\end{tabular}
\end{table}

\begin{figure}
	\centering
	\subfigure[Domain 1]{
		\centering
		\includegraphics[width=0.3\textwidth]{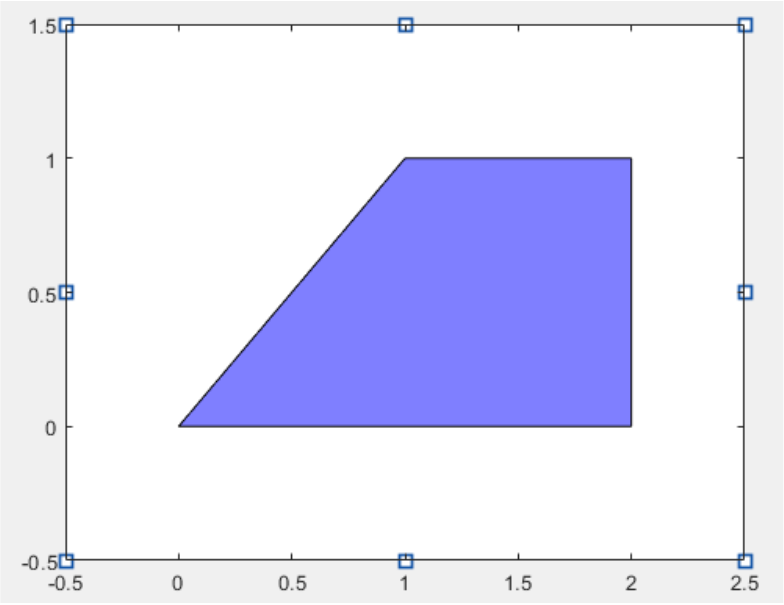}}
	\subfigure[Domain 2]{
		\centering
		\includegraphics[width=0.3\textwidth]{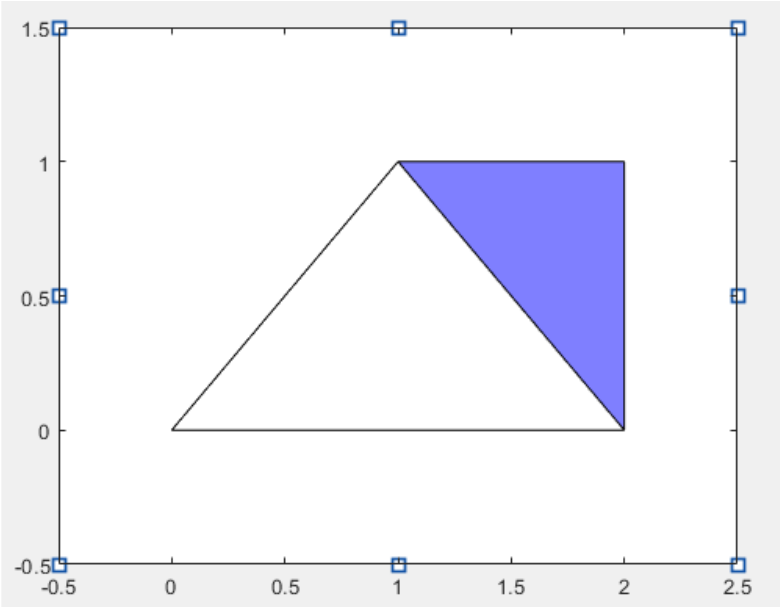}}
	\subfigure[Function]{
		\centering
		\includegraphics[width=0.3\textwidth]{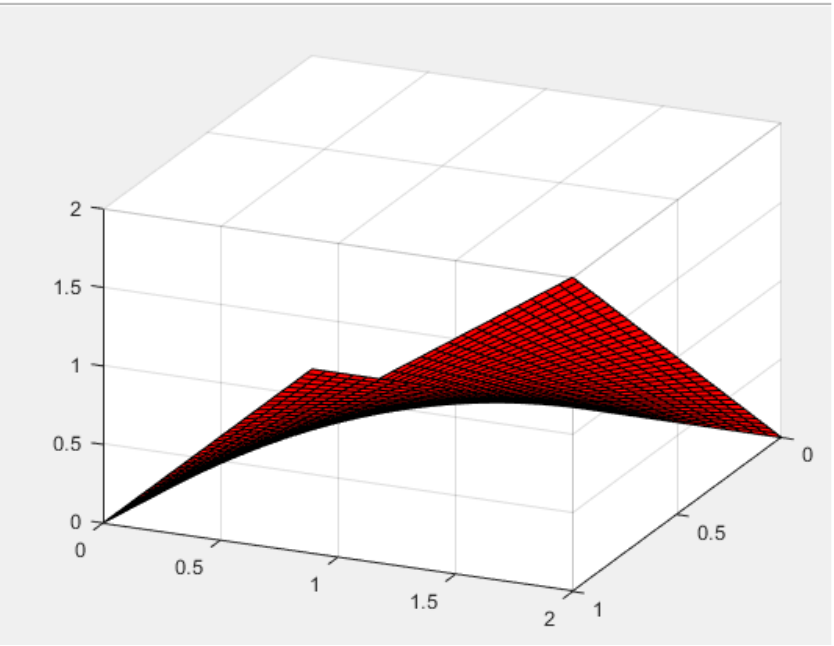}}
	\subfigure[Domain of the conjugate]{
		\centering
		\includegraphics[width=0.6\textwidth]{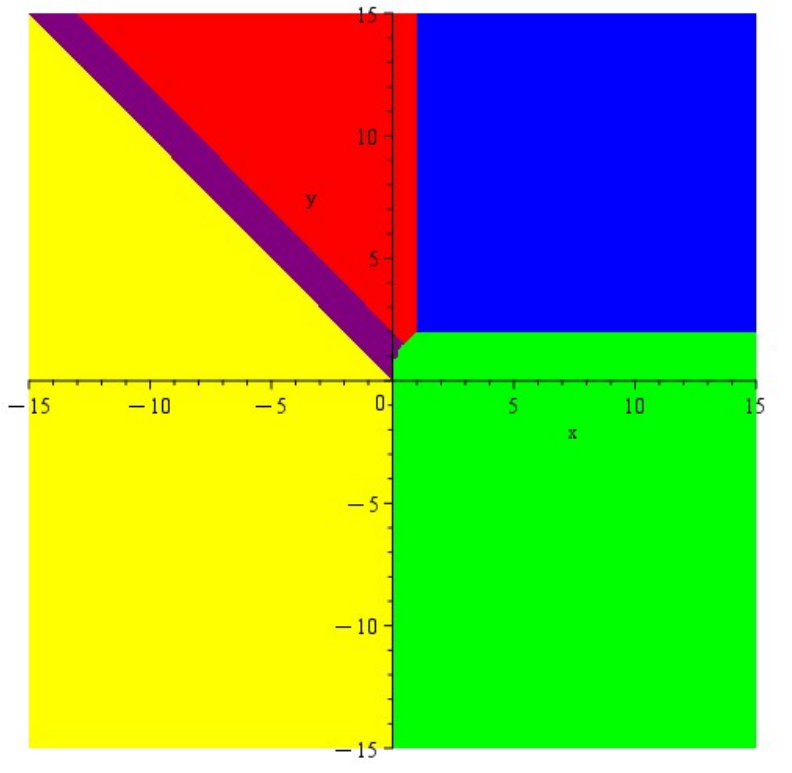}}
	\caption{Experiment 1: Computing the conjugate of a function with the domain stored as one piece gives the same result as the same function with the domain stored as two pieces.}
	\label{fig:expt1}
\end{figure}

\begin{figure}
	\centering
	\subfigure[Domain of the conjugate in all cases.]{
		\centering
		\includegraphics[width=0.35\textwidth]{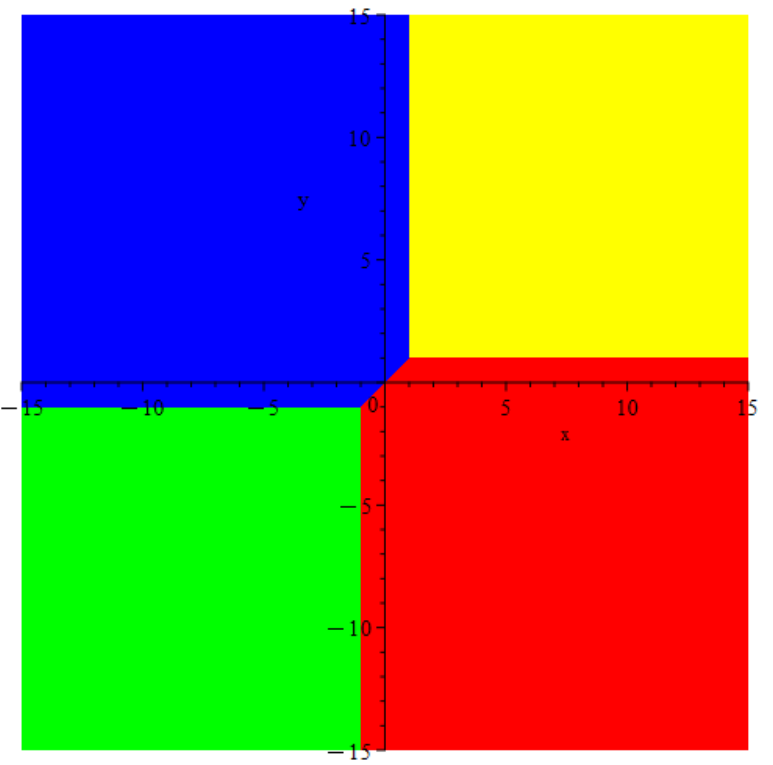}}
	\subfigure[Time to compute the conjugate vs. pieces and vs. edges.]{
		\centering
		\includegraphics[width=0.55\textwidth]{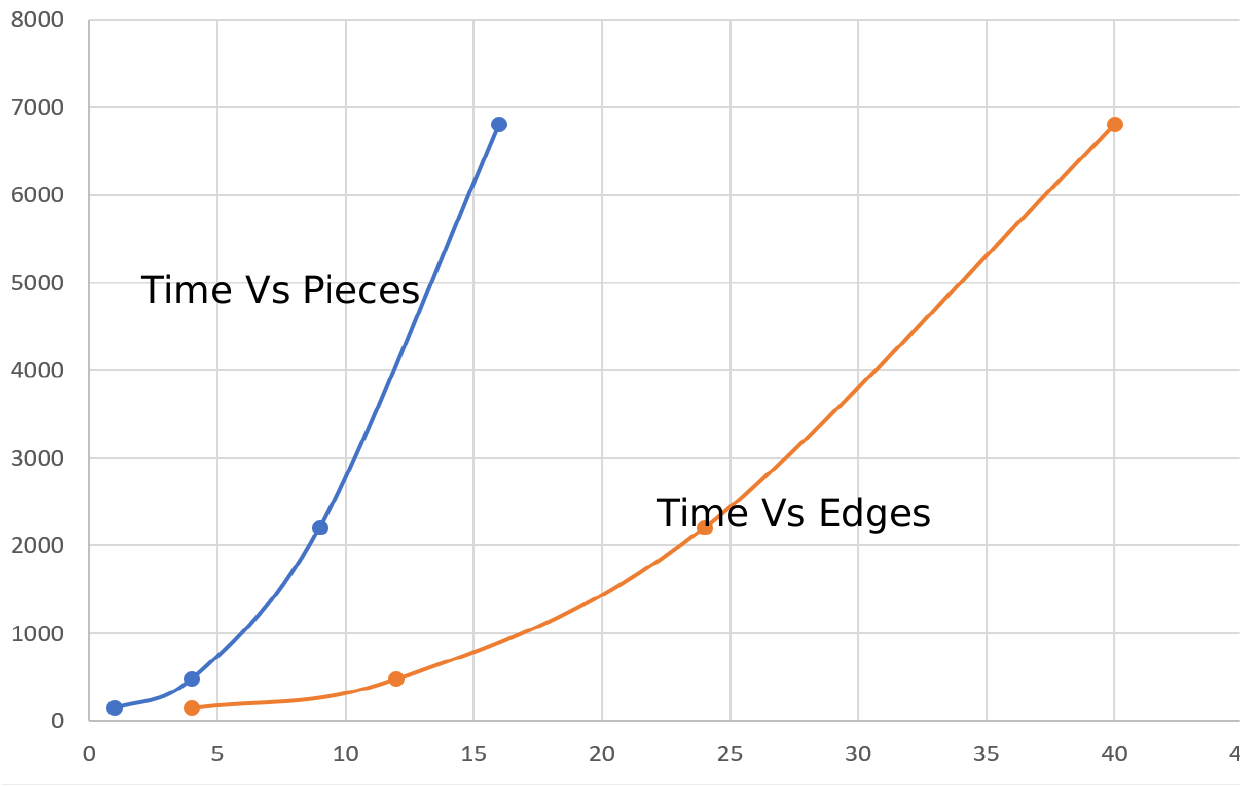}}
	\caption{Conjugate of $xy$ over unit square with each piece subdivided into $n^2$ pieces and corresponding computation time per piece and per edge.}
	\label{fig:grapg}
\end{figure}

Observing the conjugates we computed in figures \ref{fig:conjugateMax2}, \ref{fig:expt1} and \ref{fig:grapg}, we can estimate the number of pieces in the conjugate. 

\begin{conjecture}
	The number of pieces in the conjugate (computed as a maximum) is equal to the number of vertices plus the number of convex edges of the overall domain.
\end{conjecture}

\section{Examples}\label{s:examples}
We performed a lot of experiments to test and time the code. The timings are recorded on a Intel(R) Xeon(R) CPU E5-1620 v4 @ 3.50GHz 3.50 GHz machine with 64.0 GB RAM, We have run the code with MATLAB R2023B on a Windows 10 operating system. The emphasis is on correct code, not speed so the code is not optimized.

\subsection{Verification}

In order to verify the code works correctly, we first verified small problems by solving them by hand. Next we computed conjugates for examples from \cite{KUMAR-19} and verified that we obtain the same results. Finally we divided the same domain in different ways and check that we obtained the same conjugate (in the first two steps, the convex envelopes and conjugates computed are different). In all cases, we found that Step 3 computation is correctly implemented.

\subsubsection{Experiment 1}
We consider $f(x,y)=x y +I_{\co \{V_1,V_2,V_3,V_4\}}$ where $V_1 = (0,0)$, $V_2 = (1,1)$, $V_3 = (2,1)$ and $V_4 = (2,0)$. This is the example given in~\cite[Example 3.3]{KUMAR-19a}.

This function and the domain of its conjugate are shown in Figure \ref{fig:expt1}.
We repeat this experiment twice, once with one piece, second time with two pieces and obtain the same conjugate.

\subsubsection{Experiment 2}

\begin{figure}
	\centering
	\subfigure[Domain 1 ]{
		\centering
		\includegraphics[width=0.3\textwidth]{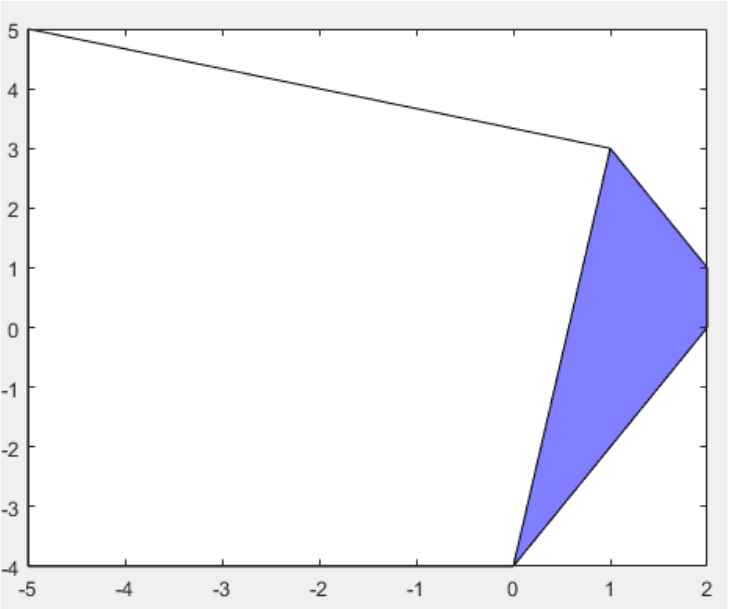}}
	\subfigure[Domain 2 ]{
		\centering
		\includegraphics[width=0.3\textwidth]{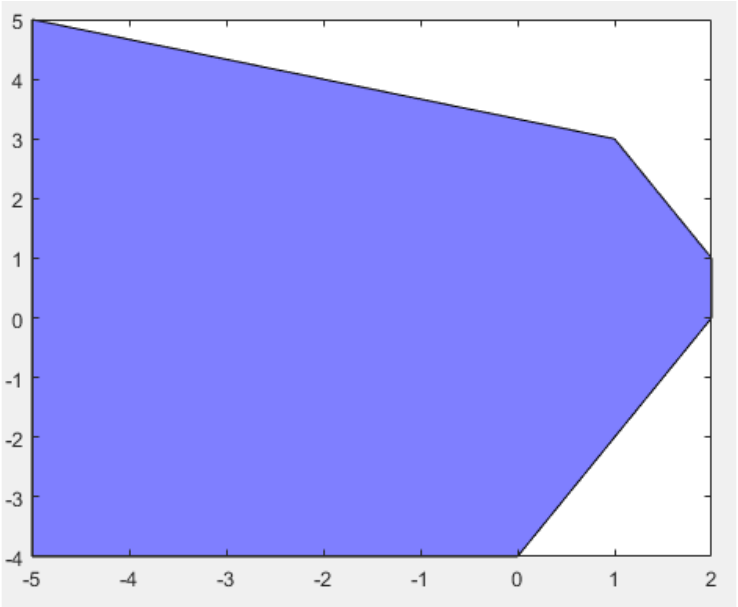}}
	
	\subfigure[Domain 3 ]{
		\centering
		\includegraphics[width=0.3\textwidth]{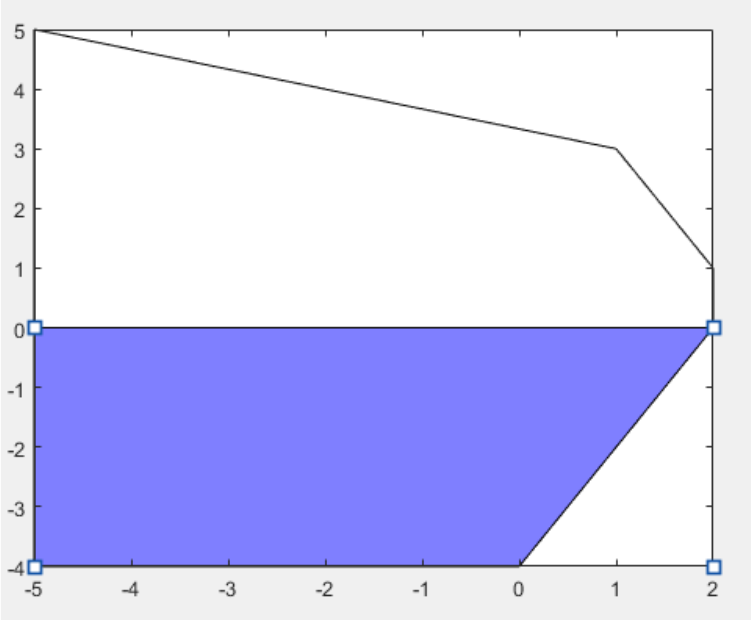}}
	\caption{Domains used in Experiment 2.}
	\label{fig:expt2}
\end{figure}

\begin{figure}
	\centering
	\subfigure[$r_{1}.$ ]{
		\centering
		\includegraphics[width=0.45\textwidth]{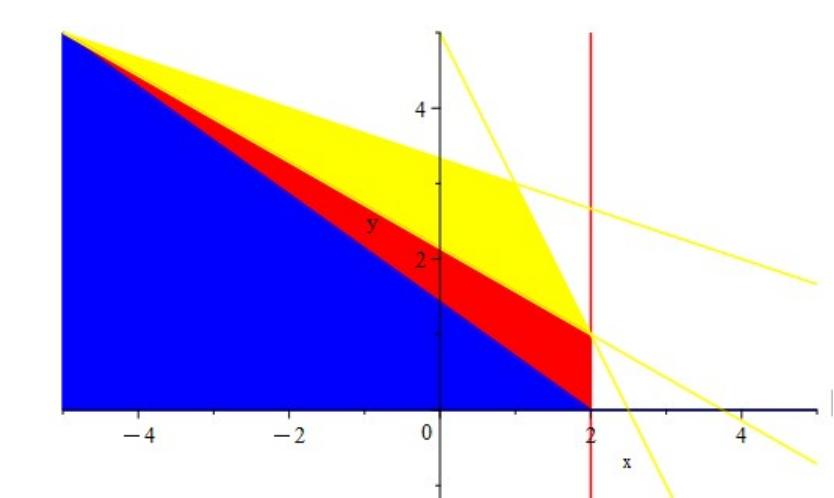}}
	\subfigure[$r_{2}.$ ]{
		\centering
		\includegraphics[width=0.45\textwidth]{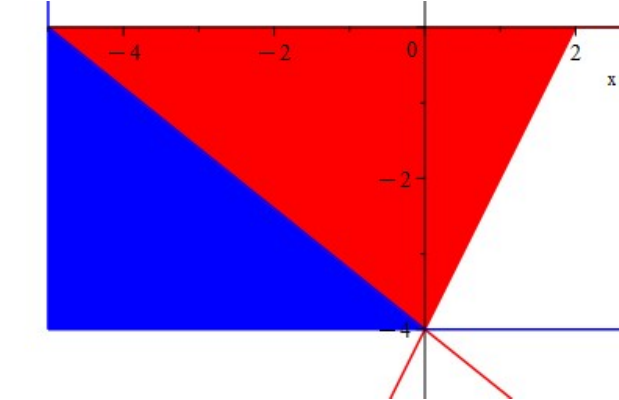}}
	\caption{Polyhedral subdivision of domain of convex envelope for Figure \ref{fig:expt2}(c).}
	\label{fig:expt2_3}
\end{figure}

We divide the same domain given in Figure~\ref{fig:example1-all} in three different ways:
\begin{enumerate}
	\item The first example has two pieces and is illustrated in Figure \ref{fig:expt2}(a). 
	\item The second example has one pieces and is illustrated in Figure \ref{fig:expt2}(b).
	\item The third example again has two pieces and is illustrated in Figure \ref{fig:expt2}(c). Polyhedral division for this example is shown in Figure \ref{fig:expt2_3}.
\end{enumerate}

Although the division is different, all three examples give the same conjugate as illustrated in Example~\ref{ex:overall2}.

\subsubsection{Timing Tests}

We find the conjugate of the function $f(x,y)=xy +I_S(x,y)$ where $S=\co \{V_1,V_2,V_3,V_4\}$ with $V_1 = (-1,-1)$, $V_2 = (-1,1)$, $V_3 = (1,1)$ and $V_4 = (-1,1)$. We then divide this piece into $n^2$ pieces and repeat the experiment where $n={1,2,3,4}$ and tabulate the results. All three tests give the same output which is given in Figure~\ref{fig:grapg}(a) and the graph of the time to find conjugate vs. number of pieces and conjugate vs edges is plotted in Figure~\ref{fig:grapg}(b).

The results for $n=1$, $n=2$ and $n=3$ are tabulated in tables \ref{table:results}, \ref{table:results2} and \ref{table:results3}. The tables contain the computation time for each step for each piece, and the average time per piece, which allow to estimate the improvement we could get by parallelisation of the code. The results are summarized in Table \ref{table:results4}. We plot graphs for time against the number of pieces originally given and the number of edges of the original domain. Both these graphs appear quadratic. 

\begin{table}[tbph]
	\centering
	\caption{Computation time to calculate the conjugate of $f(x,y)=xy +I_S(x,y)$ with $n = 1$ piece.}\label{table:results}
	\begin{tabular}{@{}lr@{}} \toprule 
        Computation step & Time (s)\\
        \midrule
		convex pieces    &	2    \\
		Time for Step 1  &	60.3 \\
		Time for Step 2  &	3.6  \\
		Time for Step 3a &	9.1  \\
		Time for Step 3b &	0.0  \\
		Total            & 	73.1 \\
		\bottomrule 
	\end{tabular}
\end{table}

\begin{table}[tbph]
	\centering
	\caption{Computation time in seconds to calculate the conjugate of $f(x,y)=xy +I_S(x,y)$ with $n = 4$ pieces; each piece is split in 2 convex pieces.
		\label{table:results2}}
	\begin{tabular}{@{}lrrrrrr@{}} \toprule 
         & \multicolumn{4}{c}{Piece No} &&\\
		 Step &\multicolumn {1}{c}{1} & \multicolumn {1}{c}{2} & \multicolumn {1}{c}{3} & \multicolumn {1}{c}{4}& Total & Avg  \\
		\toprule 
		Step 1        & 57.9	&66.5	&66.2	&55.2	&245.8  &61.5 \\
		Step 2        & 5.3	&4.9	&4.8	&2.9	&17.9	&4.5  \\
		Step 3a       & 12.6	&12.7	&12.9	&8.0	&46.2	&11.6 \\
		Average per piece &        &       &       &       &       &77.5 \\
		Step 3b       &		&       &       &       &168.4   	  \\
		Total                  & 		&       &       &       &478.3	      \\
		
		\bottomrule 
	\end{tabular}
\end{table}
\begin{table}[tbph]
	\centering
	\caption{Computation time in seconds to calculate the conjugate of $f(x,y)=xy +I_S(x,y)$ with $n = 9$ piece; each piece is split in 2 convex pieces.
		\label{table:results3}}
	\begin{tabular}{@{}lrrrrrrrrrrr@{}} \toprule 
        & \multicolumn{9}{c}{Piece No}  &&\\
		Step & \multicolumn {1}{c}1 & \multicolumn {1}{c}2 & \multicolumn {1}{c}3 & \multicolumn {1}{c}4& \multicolumn {1}{c}5 & \multicolumn {1}{c}6 & \multicolumn {1}{c}7 & \multicolumn {1}{c}8 & \multicolumn {1}{c}9 &Total & Avg  \\
		\toprule 
		Step 1        &48.4  &46.8	&48.6 &49.6	&47.8 &46.9	&46.9 &47.7	&49.1 &431.8 &47.9	\\
		Step 2        &3.4	  &3.7	&3.6  &3.5	&3.6  &3.9	&3.7  &3.5	&3.6  &32.4  &3.6	\\
		Step 3a       &9.1	  &9.9	&3.6  &9.1	&8.9  &11.2	&9.2  &9.1	&9.6  &79.5  &8.8	\\
		Average           &	  &		&	  &		&	  &		&	  &		&	  & 	 &60.4	\\
		Step 3b       &	  &		&	  &		&	  &		&	  &		&	  &1,661.4		\\
		Total                  &      &		&	  &		&	  &		&	  &		&	  &2,205.1&		\\
		
		\bottomrule 
	\end{tabular}
\end{table}


\begin{table}[tbph]\label{Timing summary in seconds.}
	\centering
	\caption{Timing summary in seconds.}\label{table:results4}
	\begin{tabular}{@{}rrrrr@{}} \toprule 
		
		nPieces	&	nVertices	&	nEdges	&	Avg time for 1 piece	&	Total time	\\ \toprule 
		1	&	4	&	4	&	73	&	146	\\
		4	&	9	&	12	&	77	&	478	\\
		9	&	16	&	32	&	60	&	2,205	\\
		16	&	25	&	40	&	58   	&	6,803	\\
		\bottomrule 
	\end{tabular}
\end{table}

Using a better algorithm, we can reduce the computation time to linear.
\begin{conjecture}
	Computing the conjugate of a PLQ function is linear with respect to the number of edges.
\end{conjecture}
Timings for both Step 1 and Step 3 can be improved significantly. In Step 1, a lot of subproblems can be discarded in advance. This code has not been implemented. In Step 3, we compute the Cartesian product of two list while computing the intersection of domains. The time taken for these loops can be reduced as this code can be implemented in linear time. Merging adjacent regions where possible would also decrease list sizes, thus giving us a faster time.

\section{Conclusion}\label{s:conclusion}

We have given a method to find the conjugate of a bivariate piecewise linear-quadratic function. We have implemented the first three steps of this method to get the conjugate of a bivariate PLQ function. It is possible to extend the code to implement Step 4 to find the biconjugate in order to get the convex envelope of the PLQ function.
In the current implementation, the complexity of Step 1 can be exponential in the worst-case, Step 2 is linear and Step 3 is polynomial. But we expect the complexity to  be reduced in the future. 

\subsection{Future Work}

The current implementation is in MATLAB. The code is written in a way that could be ported to other languages and architectures. It can be extended to incorporate the functionality that has been mentioned as future possibilities.

Future work can be pursued in several directions. Now that we have the conjugate of the PLQ function, the next natural step is to compute the biconjugate in order to obtain the convex envelope; this requires computing the conjugate of piecewise quadratic functions defined on parabolic subdivisions.

As we have implemented this computation in four steps, we can improve each of these steps. In Step 1, we can reduce the complexity of the algorithm. 	When there is no convex edge - complexity is $\binom{n}{2}$. When there is a convex edge, with a vertex opposite it - this gives a subdivision which is generally the maximum on solving the subproblems. Thus when there is a convex edge, we can reduce the number of subproblems being solved.

\begin{figure}
		\centering
		\subfigure[$r_{1,j}$ ]{
			\centering
            {\includegraphics[trim=100 220 150 220, clip,height=150px]{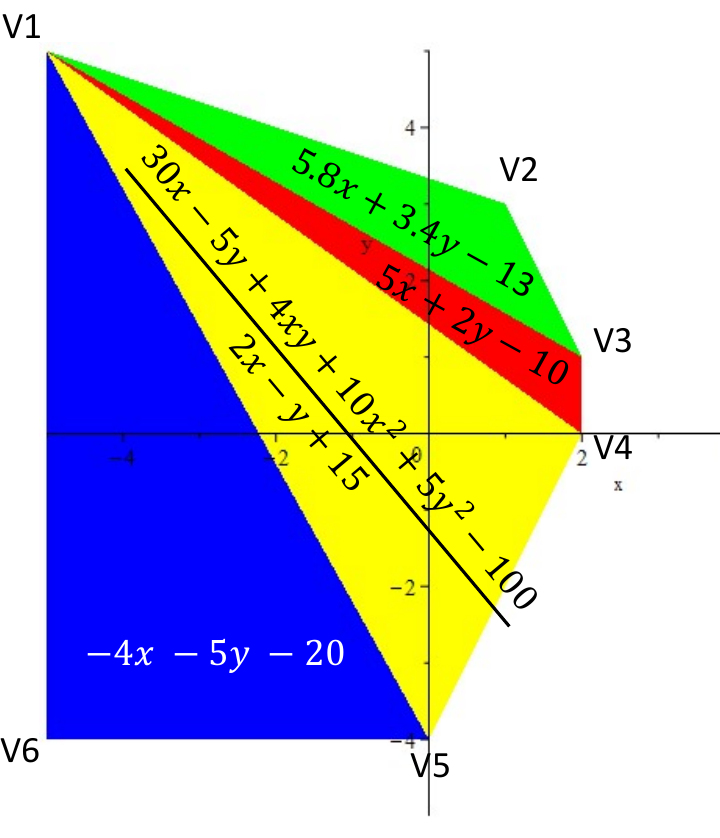}}}
		\\
		\subfigure[$s_{1,1,k}$ ]{
			\centering
			{\includegraphics[height=140px]{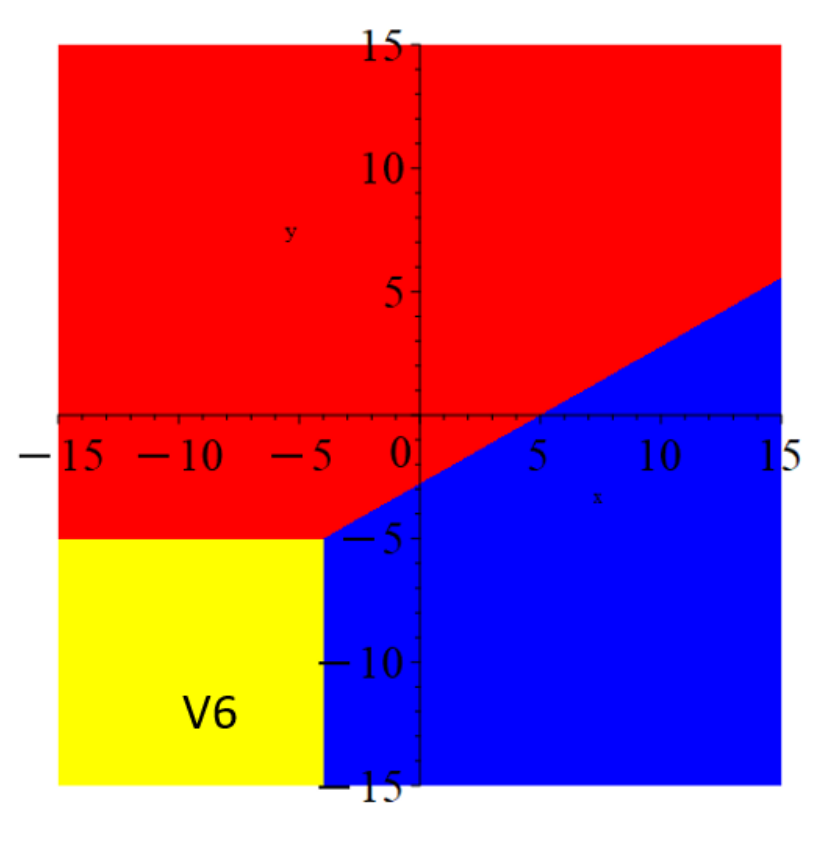}}}
		\subfigure[$s_{1,4,k}$ ]{
			\centering
            {\includegraphics[trim=50 100 150 300, clip,height=150px]{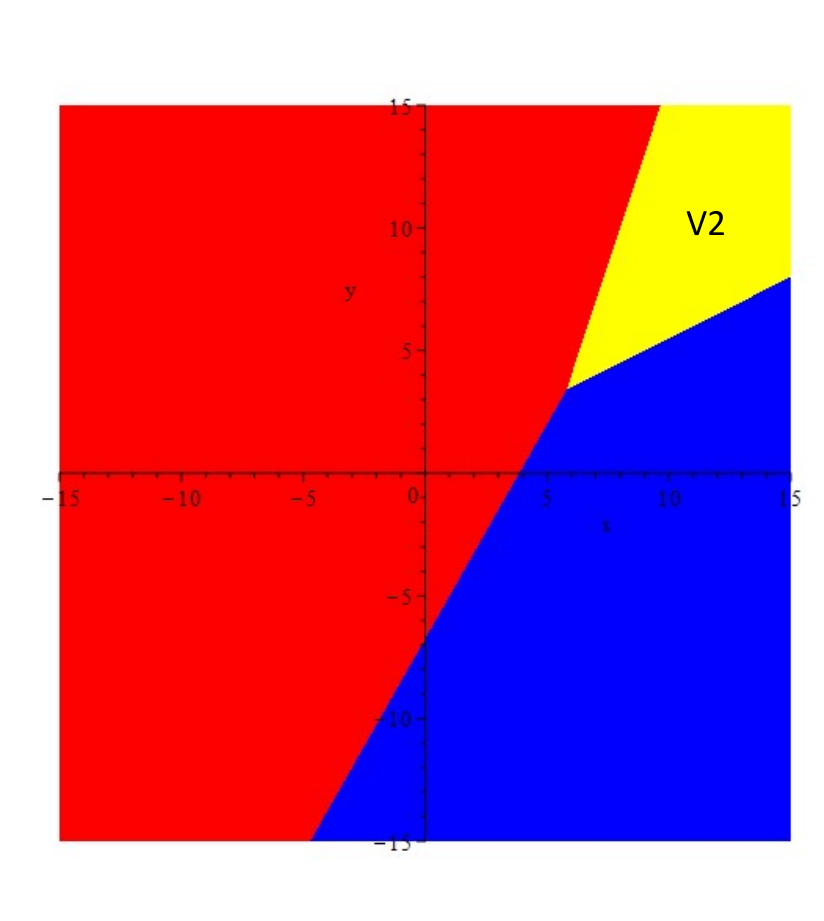}}}
		\\
		\subfigure[$s_{k}$ ]{
			\centering
            {\includegraphics[trim=180 300 120 230, clip,height=150px]{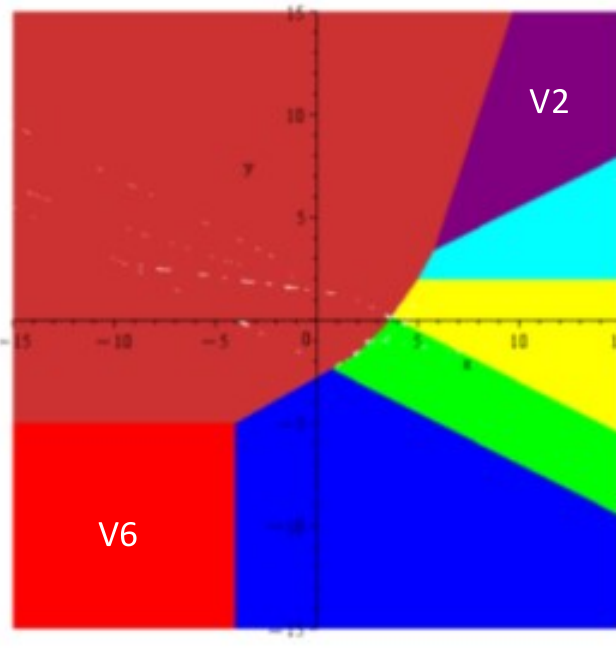}}}
		\caption{Domain of a conjugate illustrating the potential for improvement. The conjugates corresponding to vertices and edges that are not a part of the division when computing the convex envelope of each piece appear exactly the same in the final conjugate.}
		\label{fig:conjugateMax_future}
	\end{figure}

Step 3a can also be improved and implemented with a lower complexity. We have given an algorithm to implement the intersection of two domains in $O(n)$. It has also been observed that the conjugates corresponding to vertices and edges that are not a part of the division when we find the convex envelope of each piece appear exactly the same in the final conjugate. This fact can be used to develop a better algorithm to find the maximum. If you observe Figure~\ref{fig:conjugateMax_future}(a), the conjugates corresponding to the vertices ($V_2, V_6$) as illustrated in Figure~\ref{fig:conjugateMax_future}(b-c), which are not a part of the subdivision, appear exactly the same in the maximum in Figure~\ref{fig:conjugateMax_future}(d).
	
We could develop a parallel algorithm to find the convex envelope of bivariate functions over polyhedral domains. This problem being piecewise has some level of natural parallelization and we could further partition the domain depending on the hardware available. Parallelisation of the code would require more than just using parallel for statements. We would need to implement reduction and smartly distribute data in order to get an efficient algorithm. 
	
We could preprocess the original problem to get different divisions in order to solve the problem faster. The preprocessing step might vary depending on whether we use parallel methods.
	
The current examples and tests were run when the overall domain was convex. We need to run tests when the overall domain is not convex.

\backmatter

%
%
%

\bmhead{Acknowledgments}

This work was supported by Discovery grant RGPIN-2018-03928 (second author) from the Natural Sciences and Engineering Research Council of Canada (NSERC).

\section*{Declarations}

\begin{description}
	\item[Funding.] This work was supported by Discovery grant RGPIN-2018-03928 (second author) from the Natural Sciences and Engineering Research Council of Canada (NSERC).
	\item[Conflict of interest.] The authors declare that they have no Conflict of interest.
	\item[Ethics approval and consent to participate.] Not applicable.
	\item[Consent for publication.] Not applicable.
	\item[Data availability.] Not applicable.
	\item[Materials availability.] Not applicable.
	\item[Code availability.] All code and examples are available at \url{git@github.com:tanmaya11/convex.git}.
	\item[Author contribution.] Not applicable
\end{description}


%
%
%
%

\newpage
\begin{appendices}

\section{Code demonstration}\label{a:demo}
We demonstrate the code on a PLQ function with two pieces. First we create symbolic variables $x$ and $y$ and then create a function $f$ in terms of $x$ and $y$. Then we input the vertices of the polyhedral regions for the two pieces and store them in d(1) and d(2). We create two pieces p(1) and p(2), and finally the plq function.

\begin{lstlisting}
	x = sym('x')
	y = sym('y')
	f=symbolicFunction(x*y);
	d(1)=domain([-5,-4;0,-4;1,3;-5,5],x,y);
	d(2)=domain([0,-4;2,0;2,1;1,3],x,y); 
	p(1) = plq_1piece(d(1),f);
	p(2) = plq_1piece(d(2),f);
	example1 = plq(p);
\end{lstlisting}

Once we have the input, we can invoke the conjugate function to get the conjugate of the function. The input and output variables are the same as we are updating fields inside the class. 

\begin{lstlisting}
	example1 = example1.conjugate
\end{lstlisting}

Now we need the output to interpret our results.

\begin{lstlisting}
	example1.printDomainMaple;
	example1.print;
	example1.printLatex;
\end{lstlisting}

\begin{figure}[h!]
	\centering
	\includegraphics[width=12cm]{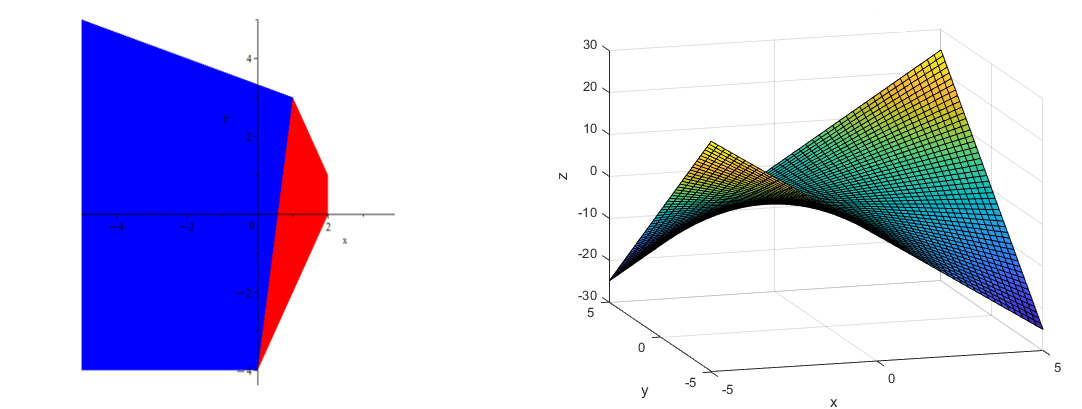}
	\caption{Function $f(x,y)=xy$ with piecewise polyhedral domain.}
	\centering
	\label{fig:10}
\end{figure}

Line 1 outputs the domain of the convex envelopes, conjugates and maximum conjugate in a format to generate diagrams in maple. The Figures created using this command are displayed in chapters (3-7).
Line 2 outputs the convex envelope, conjugate and maximum in the MATLAB command window.
Line 3 outputs the convex envelope, conjugate and maximum in \LaTeX format.
The example is plotted in Figure~\ref{fig:10}.

\section{Symbolic verification of proofs of Proposition~\ref{thm_pexpr} 
in Section~\ref{s:alg}}\label{a:proofs}

\subsection{Proof of Proposition~\ref{thm_pexpr}}\label{a:proof_them_pexpr}
The MATLAB code below verifies our computation by evaluating $\delta$ to zero. 

\begin{lstlisting}[escapeinside={(*}{*)}]
	m = sym('m')
	q = sym('q')
	
	(*$\psi_{01} = sym('\psi_{01}');$*)
	(*$\psi_{02} = sym('\psi_{02}');$*)
	(*$\psi_{03} = sym('\psi_{03}');$*)
	(*$\psi_{11} = sym('\psi_{11}');$*)
	(*$\psi_{12} = sym('\psi_{12}');$*)
	(*$\psi_{13} = sym('\psi_{13}');$*)
	(*$\psi_{21} = sym('\psi_{21}');$*)
	(*$\psi_{22} = sym('\psi_{22}');$*)
	(*$\psi_{23} = sym('\psi_{23}');$*)
	
	(*$t_0 = (-\psi_{01}-m\psi_{02})/(2(\psi_{11}+m\psi_{12}));$*)
	(*$t_1 = 1/(2(\psi_{11}+m\psi_{12}));$*)
	(*$t_2 = m/(2(\psi_{11}+m\psi_{12}));$*)
	(*$\gamma_{10} = t_1*(\psi_{23}+q\psi_{22})/(\psi_{11}+m\psi_{12});$*)
	(*$\gamma_{01} = t_2*(\psi_{23}+q\psi_{22})/(\psi_{11}+m\psi_{12});$*)
	(*$\gamma_{00} = (t_0*(\psi_{23}+q\psi_{22})-\psi_{13}-q(\psi_{12}))/(\psi_{11}+m\psi_{12});$*)
	
	(*$\zeta_{11} = -(\psi_{11}\gamma_{10}+m\psi_{12}\gamma_{10})^2/(\psi_{23}+q\psi_{22})+ \gamma_{10};$*)
	(*$\zeta_{12} = -(2(\psi_{11}\gamma_{01}+m\psi_{12}\gamma_{01})(\psi_{11}\gamma_{10}+m\psi_{12}\gamma_{10})) /(\psi_{23}+q\psi_{22}) + 
	\gamma_{01} + m\gamma_{10}$*)
	(*$\zeta_{22} = -(\psi_{11}*\gamma_{01}+m\psi_{12}\gamma_{01})^2/(\psi_{23}+q\psi_{22}) + \gamma_{01} m$*)
	
	(*$\zeta_{10} = -2(\psi_{11}\gamma_{01}+m\psi_{12}\gamma_{10})(\psi_{13}+\psi_{11}\gamma_{00}+\psi_{12}(q+m\gamma_{00}))
	/(\psi_{23}+q\psi_{22}) - m\psi_{02}\gamma_{10} + \gamma_{00} - \psi_{01}\gamma_{10};$*)
	(*$\zeta_{01} = -(2(\psi_{11}\gamma_{01}+m\psi_{12}\gamma_{01})
	(\psi_{13}+\psi_{11}\gamma_{00}+\psi_{12}(q+m\gamma_{00})))/((\psi_{23}+q\psi_{22})) 
	- m\psi_{02}\gamma_{01} - \psi_{01}\gamma_{01} + m\gamma_{00}+q;$*)
	(*$\zeta_{00} = -(\psi_{13}+\psi_{11}\gamma_{00} +\psi_{12}(q+m\gamma_{00}))^2/(\psi_{23}+q\psi_{22})
	-\psi_{03} - \psi_{01}\gamma_{00} - \psi_{02}(q+m\gamma_{00});$*)
	
	(*$D = \zeta_{12}^2 - 4\zeta_{11}\zeta_{22}$*)
	simplifyFraction(D)
	
	(*$s_1 = sym('s_1')$*)
	(*$s_2 = sym('s_2')$*)
	
	disp('quad')
	(*$fq = simplifyFraction(\zeta_{11}s_1^2 + \zeta_{12}s_1s_2 + \zeta_{22}s_2^2   + \zeta_{10}s_1 + \zeta_{01}s_2 + \zeta_{00})$*)
	
	(*$[cx,tx] = coeffs(fq,[s_1,s_2])$*)
	
	(*$\delta = simplifyFraction(cx(2)^2 - 4cx(1)cx(4))$*)
\end{lstlisting}

\subsection{Proof of zero denominator at one vertex }\label{a:nofractional}
\begin{lstlisting}[escapeinside={(*}{*)}]
% quad - linear case of step 1 
% rational function with 0/0 at vertex
a= sym('a')
b= sym('b')
m= sym('m')
q= sym('q')
fv = sym('fv')
xv = sym('x_1')
yv = sym('y_1')
fv = xv*yv
dl= sym('dl')
du= sym('du')
etah = -(a+m*b-q)^2/(4*m)-b*q
etaw = fv - a*xv - b*yv
eq = etah - etaw
z = sym('z')
% z = a+mb
eq2 = subs(eq, a, z-m*b)
bv = solve(eq2,b)
obj = etaw + a*x+ b*y
obj = subs(obj, a, z-m*b)
obj = subs(obj, b, bv)
[cz,terms] = coeffs(obj,z)
psi2 = -simplify(cz(1) )
psi1 = simplify(cz(2) )/2
psi0 = simplify(cz(3) )
% vertex 0/0 for obj1
subs(psi2, [x,y], [xv,yv])
subs(psi1, [x,y], [xv,yv])
obj1 = simplifyFraction(psi1^2/psi2 + psi0)
obj2 = simplifyFraction(-psi2*dl^2+2*dl*psi1+psi0)
obj2 = simplifyFraction(-psi2*du^2+2*du*psi1+psi0)
\end{lstlisting}
\end{appendices}


\bibliography{conjugate}

\end{document}